\documentclass[11pt,reqno]{amsart}
\usepackage[T1]{fontenc}
\usepackage[utf8]{inputenc}
\usepackage{lmodern,amsmath,amssymb,amsthm,mathtools,mathrsfs,array,booktabs,tabularx}
\usepackage[margin=28mm]{geometry}
\usepackage[colorlinks=true,linkcolor=blue,citecolor=blue,urlcolor=blue]{hyperref}
\numberwithin{equation}{section}
\newtheorem{theorem}{Theorem}[section]
\newtheorem{lemma}[theorem]{Lemma}
\newtheorem{proposition}[theorem]{Proposition}
\newtheorem{corollary}[theorem]{Corollary}
\theoremstyle{remark}\newtheorem{remark}[theorem]{Remark}
\newcommand{\R}{\mathbb R}

\newcommand{\T}{\mathbb T}
\newcommand{\one}{\mathbf1}
\newcommand{\Id}{\operatorname{Id}}
\newcommand{\diag}{\operatorname{diag}}
\newcommand{\osc}{\operatorname{osc}}
\newcommand{\spr}{\operatorname{spr}}
\newcommand{\norm}[1]{\left\lVert #1\right\rVert}
\title[Communication, profiles and critical asymptotics]{Principal profiles, characteristic communication and critical asymptotics for cooperative nonlocal systems on the circle}
\author[Bing He]{Bing He}
\address{School of Mathematics, Guangdong University of Education, Guangzhou 510303, PR China}
\email{hzs314@163.com}
\date{}
\makeatletter
\@ifundefined{subjclassname@2020}{%
  \@namedef{subjclassname@2020}{\textup{2020} Mathematics Subject Classification}}{}
\makeatother
\hypersetup{pdftitle={Principal profiles, characteristic communication and critical asymptotics for cooperative nonlocal systems on the circle},pdfauthor={Bing He}}
\begin{document}

\begin{abstract}
We study cooperative nonlocal transport on a circle with
species-dependent constant velocities and arbitrary $L^1$ probability kernels.
For periodic coefficients, the period map is a compact perturbation of
the no-jump map. We construct separate receiving and sending
characteristic scans: under a two-tree condition based at one active
component, a positive Floquet profile exists exactly when the full
period map has strictly larger spectral radius than its no-jump part.
The other jump rates may vanish. The scans provide an integral minorant
and positive jump return; a separate coefficient inequality verifies
the spectral gap. Smooth locking examples distinguish communication
from frozen-time mixing, while a resonant example separates positive
return from spectral separation. Autonomous nonzero velocities and
union-graph connectivity yield a simple principal pair; one active
component suffices for exponential mixing. We also prove local
perturbation results. For quadratic self-regulation, the same principal
projection gives a transverse response whose averaged cubic coefficient
is a nonnegative asymptotic variance. It determines both the critical
$t^{-2}\log t$ correction and the second-order period-averaged amplitude
of the positive branch.
\end{abstract}

\maketitle

\section{Introduction}
Cooperative transport and nonlocal dispersal on a compact periodic
space need not select a positive profile after the growth factor is
removed. Species-dependent velocities and inactive jump rates create
different obstructions to spatial communication. We ask when local
cooperation and nonlocal dispersal overcome these obstructions and yield
a positive principal profile with exponential mixing.

On $\T_L=\R/L\mathbb Z$, component $i$ has constant transport
coefficient $c_i\in\R$, jump rate $d_i\ge0$, and probability density
$J_i\in L^1(\R)$. A continuous Metzler matrix $B$ describes growth
and cooperative transfer. The generator acts by
\[
 (Af)_i=c_i f_i'+d_i(J_i*f_i-f_i)+\sum_jB_{ij}f_j,\qquad
 (J_i*f_i)(x)=\int_\R J_i(h)f_i(x-h)\,dh.
\]
The free evolution is $f_i(x+c_it)$; velocities may agree or have
opposite signs. Component $i$ is active when $d_i>0$. No moment,
continuity or pointwise upper bound is required of the kernels.

In the autonomous problem, nonzero velocities and strong connectivity
of the spatial union graph give a compact resolvent and a simple
positive principal pair. One active component suffices for exponential
mixing. In the periodic problem, let $M$ and $V_0$ be the full and
no-jump period maps; $V_0$ retains the diagonal losses $-d_i$.
Write $r(S)$ for the spectral radius of a bounded operator $S$.
We prove that $M-V_0$ is compact. Within either the all-active
union-graph class or the partially dispersive two-tree class,
a positive Floquet profile exists exactly when $r(M)>r(V_0)$.
Communication and positive jump return are established separately
from this spectral separation.

The two trees distinguish receiving a constant pointwise floor from
sending the mass of arbitrary input to a designated active component.
A finite convolution power in that component then redistributes the
mass around the circle. The resulting chronological integral minorant
leads, after principal normalization, to exponential mixing.

The contributions form one chain, from characteristic communication to
principal normalization and then to nonlinear relaxation:
\begin{enumerate}
\item Theorems~\ref{tp:th:main} and \ref{oa:criterion} combine a
      compact full/no-jump splitting with explicit communication
      estimates. The two-tree construction distinguishes the receiving
      and sending velocities when only one component disperses.
      Corollary~\ref{oa:row} supplies a coefficient test for the separate
      spectral gap; Proposition~\ref{prop:heterogeneous-gate} applies it
      when the principal profile is spatially nonconstant.
\item The locking and resonance results
      (Propositions~\ref{prop:smooth-lock} and \ref{prop:dual-lock},
      Theorem~\ref{th:floquet}, and
      Corollary~\ref{cor:one-active-no-gap}) distinguish the roles of
      communication, jump return and spectral separation. They explain
      why frozen autonomous mixing cannot replace a chronological
      periodic estimate.
\item Lemma~\ref{nl:lem:dissipation} and
      Theorem~\ref{nl:th:normal-form} identify the cubic amplitude
      coefficient as an asymptotic variance. The same coefficient
      controls the critical logarithmic correction and the second-order
      amplitude of the positive periodic branch. This connects the
      linear principal projection to a quantified nonlinear conclusion.
\end{enumerate}
The autonomous results provide a complementary compact-resolvent route
and allow local velocity perturbations; periodic perturbations keep the
velocities fixed.

For autonomous nonlocal operators, Berestycki--Coville--Vo \cite{BCV}
develop generalized principal values, Coville--Hamel \cite{CH} treat
scalar interval drift, and Bentout--Vo \cite[model~(1.6)]{BV} study
an interval system with a different domain. The scalar circle theory of
Ducrot--Griette--Liang \cite[Assumption~2.1, Theorem~2.3]{DGL} allows
nonvanishing $C^1$ drift; their vanishing- and sign-changing-drift
results include measure-valued eigenpairs and multiplicity
\cite[Theorems~2.9--2.10]{DGL}. Table~\ref{tab:comparison} contrasts
these hypotheses with the coupled constant-velocity setting here.

For cooperative time-periodic systems without spatial transport,
Feng--Li--Ruan--Xin \cite[Theorem~2.13, Remark~2.14]{FLRX} give
full/no-jump separation and its positive converse, including periodic
boundary conditions. Bao--Shen \cite[Theorems~2.1--2.3]{BS} prove
existence and simplicity criteria with a common regular kernel.
Wang--Zhang \cite[(L2), Lemmas~2.1--2.3]{WZ} allow distinct continuous
kernels and time-space averaged irreducibility; Wu--Li--Sun--Vo
\cite{WLSV} allow coupled integral terms. The scalar whole-space Floquet theory of
Cloez--El Abdouni--Gabriel \cite{CEAG} gives exponential $L^1$
attraction under confinement and kernel hypotheses; see
Table~\ref{tab:comparison}. The variable-drift and whole-space
settings of \cite{DGL,CEAG} extend in directions outside the present
circle problem, so the comparison is not an inclusion of their results.

\begin{table}[htbp]
\caption{Neighboring transport results and the present communication mechanism.}
\label{tab:comparison}
\small
\setlength{\tabcolsep}{4pt}
\renewcommand{\arraystretch}{1.15}
\begin{tabularx}{\textwidth}{@{}>{\raggedright\arraybackslash}p{.19\textwidth}>{\raggedright\arraybackslash}X>{\raggedright\arraybackslash}X@{}}
\toprule
Work & Setting and hypotheses & Relevant conclusion or mechanism\\
\midrule
Ducrot--Griette--Liang \cite{DGL}, Assumption 2.1 and Theorem 2.3
& Scalar autonomous circle problem; nonvanishing $C^1$ drift,
continuous potential and strictly positive continuous kernel.
& Unique normalized principal pair. Other theorems treat vanishing
and sign-changing drift, including measure-valued pairs.\\
Cloez--El Abdouni--Gabriel \cite{CEAG}, Theorem 1.1
& Scalar time-periodic problem on $\R^d$ in $L^1$; confining selection,
controlled characteristic residence, kernel nonconcentration and positivity.
& Floquet pair and exponential attraction using positive-operator
spectral methods.\\
Present work, Theorems~\ref{tp:th:main} and \ref{oa:criterion}
& Coupled circle system on $C$; species-dependent constant velocities,
arbitrary $L^1$ convolution densities, possibly inactive components.
& Compact full/no-jump splitting; receiving and sending scans;
exact radius criterion within the stated communication classes.\\
Present work, Lemma~\ref{nl:lem:dissipation} and
Theorem~\ref{nl:th:normal-form}
& Quadratic self-regulation after principal normalization and
exponential mixing.
& Cubic coefficient identified as a stationary asymptotic variance;
critical and positive-branch second-order expansions.\\
\bottomrule
\end{tabularx}
\end{table}

The operator argument combines established evolution-semigroup,
positive-inverse and perturbation theory \cite[Propositions~3.1--3.2]{HR},
\cite[Theorems~3.1, 4.7]{GM}, \cite{EN,Kato}, with principal
normalization and Doeblin contraction \cite{HM,CG,BCGM}.
The model-specific contributions are the two characteristic scans,
their locking and resonance obstructions, and the identification of
the nonlinear cubic coefficient with an asymptotic variance.

The paper is organized as follows. Section~\ref{sec:main} states the
autonomous results and summarizes the periodic and nonlinear results.
Section~\ref{sec:autonomous-theory} proves the autonomous principal,
mixing and perturbation theorems. Section~\ref{sec:periodic-theory}
establishes the abstract compact-splitting criterion, its periodic
transport realization and the two-tree estimates.
Section~\ref{sec:application} develops quadratic self-regulation,
including the variance identity and second-order expansions.
Section~\ref{sec:examples-scope} collects the communication and drift
examples. Appendix~\ref{app:technical} gives the spectral obstructions
and a direct proof of the periodic realization;
Appendix~\ref{app:applications} contains the remaining nonlinear and
drift calculations.

\section{Framework and main results}\label{sec:main}

\subsection{Autonomous setting and principal data}

Let $L>0$, $m\ge1$, and $\T_L=\R/L\mathbb Z$, equipped with Lebesgue
measure of total mass $L$. On
\[
 X=C(\T_L;\R^m),\qquad
 \norm f=\max_i\sup_{x\in\T_L}|f_i(x)|,
\]
consider
\begin{equation}\label{eq:A}
 (Af)_i=c_i f_i'+d_i(K_if_i-f_i)+\sum_{j=1}^mB_{ij}(x)f_j,
 \qquad
 K_i u(x)=\int_\R J_i(h)u(x-h)\,dh.
\end{equation}
The autonomous assumptions used here and in Section~\ref{sec:autonomous-theory} are:
\begin{enumerate}
\item Every $c_i$ is a nonzero real constant, and $d_i\ge0$.
\item $J_i\in L^1(\R)$, $J_i\ge0$ a.e., and $\int_\R J_i=1$.
\item $B\in C(\T_L;\R^{m\times m})$ is Metzler:
      $B_{ij}(x)\ge0$ for $i\ne j$.
\item The directed graph with edge $j\to i$ whenever
      $i\ne j$ and $I_{ij}:=\int_0^LB_{ij}(x)\,dx>0$ is strongly connected.
      For $m=1$ this condition is vacuous.
\end{enumerate}
Order on $X$ is componentwise. We write $f\gg0$ when
$\min_{i,x}f_i(x)>0$, $f\ll0$ when $-f\gg0$, and $\one$ for the
constant vector of ones. The notation $\diag$ forms a diagonal matrix.
For Banach spaces $E,F$, $\mathcal L(E,F)$ denotes the bounded linear
operators from $E$ to $F$.
The velocities need not be equal or have the same sign. No boundedness,
continuity, support or moment condition on a density $J_i$ is imposed.
All coefficients are autonomous. Write
\begin{equation*}
 \begin{gathered}
 D(A)=D=C^1(\T_L;\R^m),\qquad
 C=\diag(c_i),\quad D_d=\diag(d_i),\\
 c_-:=\min_i|c_i|>0,\quad c_+:=\max_i|c_i|.
 \end{gathered}
\end{equation*}
We use the complexification when taking spectra. For an operator $H$,
$\sigma(H)$ and $\rho(H)$ denote its spectrum and resolvent set;
$s(H)=\sup\{\operatorname{Re}z:z\in\sigma(H)\}$ is its spectral bound.
For a bounded operator, $r(H)$ (also written $\spr(H)$) is its spectral
radius. The $C^1$ norm is $\norm f_{C^1}=\norm f+\norm{f'}$.
For positive scalars, $f\asymp g$ means that $f/g$ is bounded above
and away from zero. Write $K=\diag(K_i)$ for the componentwise convolution operator and
$\Id$ for the identity.

The formal adjoint is the separate realization
\begin{equation}\label{eq:adjoint}
 A^\dagger=-C\partial_x+D_d(K^{\mathrm{ref}}-\Id)+B^T,
 \qquad D(A^\dagger)=D,
 \qquad J_i^{\mathrm{ref}}(h)=J_i(-h).
\end{equation}
Here $K^{\mathrm{ref}}$ is the componentwise convolution with the
reflected kernels $J_i^{\mathrm{ref}}$, and $B^T$ is the transpose matrix.
This realization is distinct from the Banach dual of $X$, whose elements
are measures. For normalization use the fixed spatial and component average
\[
 \mathcal N(f)=\frac1{mL}\sum_i\int_0^Lf_i(x)\,dx,
 \qquad \mathcal P=\{z\in D:\min_{i,x}z_i(x)>0\}.
\]

For a semigroup $S$, its growth bound $\omega_0(S)$ is the infimum of
$\omega\in\R$ for which $\|S(t)\|\le C_\omega e^{\omega t}$ for all
$t\ge0$ and some $C_\omega<\infty$.

\begin{theorem}[Principal pair and comparison]\label{th:principal}
Under the stated assumptions, $A$ generates a positive $C_0$ semigroup
$S(t)$ and has compact resolvent. There are $s\in\R$ and
$\phi,\psi\in D$, with every component strictly positive, such that
\begin{equation}\label{eq:eigen}
 A\phi=s\phi,\qquad A^\dagger\psi=s\psi,\qquad
 \mathcal N(\phi)=1,\qquad \int_0^L\psi^T\phi\,dx=1.
\end{equation}
The eigenvalue $s$ is isolated and algebraically simple. It equals the
spectral bound and the semigroup growth bound. If
$\phi_{\min}=\min_{i,x}\phi_i(x)$ and
$\phi_{\max}=\max_{i,x}\phi_i(x)$, then
\begin{equation}\label{eq:norm-pair}
 e^{st}\le\norm{S(t)}\le\frac{\phi_{\max}}{\phi_{\min}}e^{st},
 \qquad
 \int_0^L\psi^TS(t)f\,dx=e^{st}\int_0^L\psi^Tf\,dx.
\end{equation}
Moreover,
\begin{equation}\label{eq:CW}
 s=\sup_{z\in\mathcal P}\min_{i,x}\frac{(Az)_i(x)}{z_i(x)}
  =\inf_{z\in\mathcal P}\max_{i,x}\frac{(Az)_i(x)}{z_i(x)}.
\end{equation}
Both extrema are attained at $\phi$. The following conditions are equivalent:
\begin{equation}\label{eq:MP}
 s<0;\qquad
 [f\in D,\ Af\ge0\Longrightarrow f\le0];\qquad
 0\in\rho(A)\ \text{and }(-A)^{-1}\ge0.
\end{equation}
\end{theorem}

The maximum principle in \eqref{eq:MP} uses the periodic generator
domain. The operator $A^\dagger$ in \eqref{eq:eigen} is the separate
$C^1$ realization \eqref{eq:adjoint}; the Banach dual acts on measures.

\begin{theorem}[Exponential principal projection]\label{th:mix}
If in addition $d_a>0$ for some species $a$, then there exist
$C_{\mathrm{mix}}<\infty$ and $\kappa>0$ such that
\begin{equation}\label{eq:mix}
 \left\|e^{-st}S(t)f-\phi\int_0^L\psi^Tf\,dx\right\|
       \le C_{\mathrm{mix}}e^{-\kappa t}\norm f\qquad(t\ge0).
\end{equation}
The raw minorant is obtained from finite graph paths, drift magnitudes,
coefficient bounds and a finite convolution minorant for $J_a$.
The normalized constants also depend on an upper bound for $s$ and on
$\phi_{\max}/\phi_{\min}$, as specified in
\eqref{eq:mix-constants}; these bounds are locally uniform under
Theorem~\ref{th:perturb}.
\end{theorem}

An active jump suffices for mixing. Theorem~\ref{th:principal}
also allows all $d_i=0$, a class in which principal functions can exist
without normalized convergence.

For perturbations, fix $m,L$ and an operator satisfying the above assumptions.
An admissible perturbed operator has the same form, with nonnegative rates,
probability densities and a continuous Metzler matrix. Use
\begin{equation}\label{eq:Delta}
 \Delta=\max_i|\widetilde c_i-c_i|
 +\max_i|\widetilde d_i-d_i|
 +\max_{ij}\norm{\widetilde B_{ij}-B_{ij}}_\infty
 +\max_i\norm{\widetilde J_i-J_i}_{L^1(\R)}.
\end{equation}

\begin{theorem}[Local perturbations and first variations]\label{th:perturb}
There are $\Delta_0>0$ and $C_{\mathrm{par}}<\infty$, depending on the fixed realization, such that every
admissible $\Delta<\Delta_0$ has nonzero component velocities and a strongly
connected union graph. Its normalized principal data satisfy
\begin{equation}\label{eq:parameter}
 |\widetilde s-s|+\norm{\widetilde\phi-\phi}_{C^1}
                 +\norm{\widetilde\psi-\psi}_{C^1}
       \le C_{\mathrm{par}}\Delta.
\end{equation}
If $d_a>0$ at the fixed realization, the constants in \eqref{eq:mix} can
be chosen common to a sufficiently small admissible neighborhood.
All these local estimates are uniform over the simultaneous spatial phases
$B(x+\vartheta)$, $\vartheta\in\T_L$, where $\vartheta$ is a spatial
phase and each perturbation is measured relative to that phase.

For an admissible kernel path, $J_i\ge0$, $\int J_i=1$ and
$\int\dot J_i=0$; the derivative below is in the ambient $L^1$
space, restricted to such paths. The normalized principal data are
locally $C^1$ in the data norms in
\eqref{eq:Delta}, in the sense of restriction of a $C^1$ branch on an
open neighborhood in the ambient real Banach space. For a differentiable
admissible parameter path, where $B$ is differentiated
in the uniform norm and $J_i$ in $L^1$, the principal value is differentiable
and
\begin{equation}\label{eq:HF}
 \begin{gathered}
 \dot s=\sum_i\int_0^L\psi_i
 \left[\dot c_i\phi_i'
  +\dot d_i(K_i\phi_i-\phi_i)
  +d_i\dot K_i\phi_i+\sum_j\dot B_{ij}\phi_j\right]dx,\\
 (\dot K_i u)(x)=\int_\R\dot J_i(h)u(x-h)\,dh.
 \end{gathered}
\end{equation}
\end{theorem}

The neighborhood is local to a fixed realization and stays away from
$c_i=0$, zero required graph integrals and a zero active rate.

\subsection{Periodic communication and nonlinear asymptotics: an overview}
Theorems~\ref{tp:th:main} and \ref{oa:criterion} establish
\[
 \text{positive Floquet profile}\quad\Longleftrightarrow\quad
 r(M)>r(V_0)
\]
under, respectively, all-active union-graph communication or the
two-tree condition based at one active component. Velocities may
vanish and the remaining rates may be zero. Both classes have an
integral minorant; Theorem~\ref{thm:compact-return} then gives
principal mixing under separation. Corollary~\ref{oa:row} verifies
the gap from coefficients.

For quadratic self-regulation, Theorem~\ref{nl:th:application} gives
threshold dynamics. At zero principal exponent, the phasewise centered
response of Lemma~\ref{nl:lem:response-domain} identifies the averaged
cubic coefficient as a nonnegative asymptotic variance
(Lemma~\ref{nl:lem:dissipation}). Theorem~\ref{nl:th:normal-form}
then gives the critical $t^{-2}\log t$ correction and the second-order
positive periodic branch. Both use the same principal projection.

\section{Autonomous principal theory and perturbations}\label{sec:autonomous-theory}

This section proves the autonomous principal-pair and mixing results,
then establishes the perturbation estimates used in the periodic analysis.

\subsection{Compact resolvents and principal functions}\label{sec:principal}

Each $K_i$ is a positive contraction on $C(\T_L)$ and preserves constants.
It maps continuous functions to continuous functions because
\[
 |K_i f(x+h)-K_i f(x)|\le\sup_y|f(y+h)-f(y)|.
\]
The diagonal translations $(T(t)f)_i(x)=f_i(x+c_i t)$ form a strongly
continuous group. Their generator domain is exactly $D$: a norm difference
quotient gives the continuous point derivative of each component, and the
converse follows from the fundamental theorem of calculus and uniform
continuity of $f_i'$. Every other term in \eqref{eq:A} is bounded on $X$.
The bounded perturbation theorem \cite{EN} therefore gives a semigroup
with domain $D$.

Choose the scalar compensation constant $H_0$ and put
\begin{equation*}
 \begin{gathered}
 H_0=1+\max_i(d_i+\norm{B_{ii}}_\infty),\quad
 H=C\partial_x-H_0\Id,\\
 Q=D_dK+B-D_d+H_0\Id\ge0,\quad
 b=\max_{i,x}\sum_jB_{ij}(x).
 \end{gathered}
\end{equation*}
The Dyson series relative to $H$ is positive term by term. In addition,
$A\one=B\one\le b\one$, so positive comparison gives
\begin{equation*}
 S(t)\one\le e^{bt}\one,\qquad
 \norm{S(t)}=\norm{S(t)\one}\le e^{bt}.
\end{equation*}
For real $\lambda>\max\{b,\norm Q-H_0,0\}$ let
$R_\lambda=(\lambda-A)^{-1}$ and $R_\lambda^0=(\lambda-H)^{-1}$.
The Laplace formula and the positive resolvent expansion give
\begin{equation}\label{eq:resolvent}
 \begin{split}
 R_\lambda&=\int_0^\infty e^{-\lambda t}S(t)\,dt
       =\sum_{n=0}^\infty(R_\lambda^0Q)^nR_\lambda^0,\\
 (R_\lambda^0f)_i(x)&=\int_0^\infty
            e^{-(\lambda+H_0)t}f_i(x+c_i t)\,dt,
 \qquad \norm{R_\lambda}\le(\lambda-b)^{-1}.
 \end{split}
\end{equation}
The series converges since $\norm{R_\lambda^0Q}<1$.
The resolvent equation for $u=R_\lambda f$ is
\begin{equation*}
 c_i u_i'=\lambda u_i-d_i(K_i u_i-u_i)-\sum_jB_{ij}u_j-f_i.
\end{equation*}
Hence $R_\lambda:X\to C^1$ is bounded, and its action into $X$ is compact.
The transport semigroup itself need not be compact.

\begin{lemma}[Scalar resolvent minorant]\label{le:scalar}
For $P_i=L/|c_i|$ and
$a_i=e^{-(\lambda+H_0)P_i}/|c_i|>0$,
\begin{equation*}
 (R_\lambda^0f)_i(x)\ge a_i\int_0^Lf_i(y)\,dy\qquad(f\ge0).
\end{equation*}
\end{lemma}
\begin{proof}
Retain only $0\le t\le P_i$ in \eqref{eq:resolvent}, replace its exponential
by the value at $P_i$, and integrate over one complete circle. A negative
velocity reverses orientation and has the same Jacobian $|c_i|$.
\end{proof}

\begin{lemma}\label{le:strong-res}
There is $\rho>0$ such that
\begin{equation*}
 (R_\lambda f)_i(x)\ge\rho\sum_j\int_0^Lf_j(y)\,dy
                 \qquad(f\ge0,\ i,x).
\end{equation*}
\end{lemma}
\begin{proof}
For each ordered pair $(i,j)$ choose a simple path
$j=i_0\to i_1\to\cdots\to i_l=i$, with $l\le m-1$, using the
empty path when $i=j$. In the $l$th term of \eqref{eq:resolvent}, retain
only its off-diagonal multiplications along this path. Apply
Lemma~\ref{le:scalar} to the input component and after every edge. Its
contribution is at least
\[
 \rho_{ij}\int_0^Lf_j,
 \qquad
 \rho_{ij}=a_*^{l+1}\prod_{r=1}^lI_{i_ri_{r-1}}>0,
 \qquad a_*=\min_i a_i.
\]
These are bounds on separate input-component entries of a positive operator,
so they can be summed over $j$. Take $\rho=\min_{i,j}\rho_{ij}$.
\end{proof}

Every nonzero continuous nonnegative vector has positive total integral.
Thus Lemma~\ref{le:strong-res} makes $R_\lambda$ strongly positive.
Moreover $R_\lambda\one\ge\rho mL\one$, so its spectral radius is positive
by iteration and the spectral-radius formula. The compact strongly positive
Krein--Rutman theorem gives an algebraically simple eigenvalue
$r=\spr(R_\lambda)>0$ and $\phi\gg0$ with $R_\lambda\phi=r\phi$.
It follows that $\phi\in D$ and
\[
 A\phi=s\phi,\qquad s=\lambda-r^{-1}.
\]
The nonzero-resolvent spectral correspondence preserves isolation and
algebraic multiplicity. Normalize $\mathcal N(\phi)=1$.
Positivity and the eigenfunction identity imply
\[
 |S(t)f|\le\frac{\norm f}{\phi_{\min}}S(t)\phi
       =\frac{\norm f}{\phi_{\min}}e^{st}\phi.
\]
Testing on $\phi$ gives the reverse norm bound in \eqref{eq:norm-pair}.
Consequently the growth bound is $s$. The Laplace resolvent for
$\operatorname{Re}z>s$ excludes spectrum to its right, whereas $s$ is
an eigenvalue. Hence the spectral bound is also $s$.

The separate operator \eqref{eq:adjoint} satisfies all the same assumptions:
its velocities are $-c_i\ne0$, its kernels are reflected probability
densities, and its graph is transposed. The preceding construction gives
$\chi\gg0$ and $s^\dagger\in\R$ with $A^\dagger\chi=s^\dagger\chi$.
For $f,g\in D$, periodic integration by parts and Fubini give
\begin{equation}\label{eq:pairing}
 \int_0^Lg^TAf\,dx=\int_0^L(A^\dagger g)^Tf\,dx.
\end{equation}
In particular $(s-s^\dagger)\int\chi^T\phi=0$. Thus $s^\dagger=s$;
set $\psi=\chi/\int\chi^T\phi$. For $f\in D$, differentiating
$\int\psi^TS(t)f$ proves the pairing in \eqref{eq:norm-pair}; density
extends it to $f\in X$.

For $z\in\mathcal P$, identity \eqref{eq:pairing} makes $s$ the weighted
average of $(Az)_i/z_i$ with strictly positive weights $\psi_i z_i$.
Its minimum is therefore at most $s$ and its maximum at least $s$.
Choosing $z=\phi$ proves both attained formulas in \eqref{eq:CW}.

If $s<0$, the integral $R_*f=\int_0^\infty S(t)f\,dt$ converges in norm
and is the positive inverse of $-A$. Indeed the finite-time generator
identity and closedness give $-AR_*f=f$, and the domain identity gives
$R_*(-Af)=f$ for $f\in D$. If $Af\ge0$, then
\[
 S(t)f-f=\int_0^tS(r)Af\,dr\ge0;
\]
norm decay implies $f\le0$. If $s\ge0$, the positive $\phi$ violates
this maximum principle. A positive $(-A)^{-1}$ is impossible when $s=0$;
when $s>0$ it would map $\phi$ to $-\phi/s$. This proves
Theorem~\ref{th:principal}.

\subsection{A fixed-time minorant and exponential convergence}\label{sec:mix}

Let
\[
 S_0(t)=e^{t(C\partial_x+B-D_d)},\qquad
 S_0^\dagger(t)=e^{t(-C\partial_x+B^T-D_d)}.
\]
These denote the generated no-jump semigroups, with the multiplication
matrix retained; no commutation with that matrix is asserted.
Both are positive, and $S(t)\ge S_0(t)$. For $m>1$ put
\begin{equation*}
 \eta=\min\left\{1,\min_{I_{ij}>0,\ i\ne j}\frac{I_{ij}}{c_+}\right\}>0,
 \qquad t_0=\frac{mL}{c_-},\qquad k_0=e^{-H_0t_0}\eta^{m-1}.
\end{equation*}
For $m=1$ set $\eta=1$ and use the same $t_0,k_0$. Write $e_j$ for the
constant vector supported in component $j$ with value one there.

\begin{lemma}\label{le:scan}
For $t\ge t_0$,
\begin{equation*}
 (S_0(t)e_j)_i(x)\ge e^{-H_0t}\eta^{m-1},\qquad
 (S_0^\dagger(t)e_j)_i(x)\ge e^{-H_0t}\eta^{m-1}.
\end{equation*}
For any fixed species $a$ and $f\ge0$ it follows that
\begin{equation*}
 S_0(t_0)e_a\ge k_0\one,\qquad
 \int_0^L(S_0(t_0)f)_a\,dx\ge k_0\sum_j\int_0^Lf_j\,dx.
\end{equation*}
\end{lemma}
\begin{proof}
Variation of constants relative to $C\partial_x-H_0\Id$, dropping all
nonnegative terms but an edge $k\to i$, gives
\begin{equation*}
 (S_0(t)f)_i(x)\ge\int_0^t e^{-H_0\tau}B_{ik}(x+c_i\tau)
                   (S_0(t-\tau)f)_k(x+c_i\tau)\,d\tau.
\end{equation*}
The no-switch contribution gives $(S_0(t)e_j)_j\ge e^{-H_0t}$.
For a simple path $j=i_0\to\cdots\to i_l=i$, let
$H=\sum_{r=1}^lP_{i_r}$, $P_i=L/|c_i|$. Induction proves, whenever
$t\ge H$,
\[
 (S_0(t)e_j)_i(x)\ge e^{-H_0t}
                     \prod_{r=1}^l\frac{I_{i_ri_{r-1}}}{|c_{i_r}|}.
\]
For the last edge, restrict the estimate above to $0\le\tau\le P_i$.
The remaining time suffices for the prefix bound, which is uniform over
its spatial argument. The last integral is exactly
$I_{ik}/|c_i|$. This uses the destination velocity, with its absolute
value; opposite signs merely reverse a complete cycle. Each factor is
at least $\eta$, $l\le m-1$ and $H<t_0$, proving the first displayed inequality. The same induction on the transposed graph
with velocities $-c_i$ proves the second.

For domain data, differentiation at complementary times and periodic
integration by parts give
\[
 \int g^TS_0(t)f=\int(S_0^\dagger(t)g)^Tf.
\]
Density extends it to all continuous data. Taking $g=e_a$ gives the displayed integral assertion. The other is already proved.
\end{proof}

The lemma scans constant species inputs, not arbitrary vector data
along a common translation. The uniform prefix bound is established
before the next edge is scanned; the edge supports may be disjoint.

\begin{lemma}\label{le:convolution}
Every nonnegative $L^1(\R)$ probability density has a periodic convolution
power satisfying
\begin{equation*}
 K^N f(x)\ge\epsilon_J\int_0^Lf(y)\,dy\qquad(f\ge0)
\end{equation*}
for some integer $N\ge1$ and $\epsilon_J>0$.
\end{lemma}
\begin{proof}
The periodization $j(x)=\sum_{k\in\mathbb Z}J(x+kL)$ is an $L^1(\T_L)$
probability density. There are $a>0$ and a measurable set $E$ of positive
measure with $j\ge a1_E$ a.e. The convolution $1_E*1_E$ is continuous
by $L^2$ translation continuity and Cauchy--Schwarz; it has positive
integral. Thus $j^{*2}\ge b1_I$ on an arc $I$ of length $\ell\in(0,L)$,
for some $b>0$. Choose $n$ with $n\ell>L$. The $n$-fold real convolution
of an interval lift of $I$ is continuous and positive on the interior
of its sum interval. A closed interval of length $L$ lies strictly in
that interior, so its periodization has a positive minimum. Consequently
$j^{*(2n)}\ge\epsilon_J>0$ a.e. Taking $N=2n$ and convolving gives the asserted minorant at every $x$.
\end{proof}

Fix $a$ with $d_a>0$, and take $N,\epsilon_J$ from
Lemma~\ref{le:convolution} for $J_a$. The diagonal generator
$G=C\partial_x-H_0\Id+D_dK$ satisfies $A-G=B-D_d+H_0\Id\ge0$.
For $\tau>0$, its $a$th transport commutes with $K_a$. Retaining
the $N$-jump term gives
\begin{equation}\label{eq:active}
 (e^{\tau G}f)_a(x)\ge e^{-H_0\tau}\frac{(d_a\tau)^N}{N!}
       (K_a^Nf_a)(x+c_a\tau)
 \ge\gamma_\tau\int_0^Lf_a,
 \qquad \gamma_\tau=e^{-H_0\tau}\frac{(d_a\tau)^N}{N!}\epsilon_J.
\end{equation}
The bounded perturbations $D_dK\ge0$ and $A-G\ge0$ have
positive Dyson terms. Their zeroth terms give $S(t_0)\ge S_0(t_0)$
and $S(\tau)\ge e^{\tau G}$. On the consecutive blocks
$t_0,\tau,t_0$, positivity and the semigroup identity give
\[
 S(2t_0+\tau)=S(t_0)S(\tau)S(t_0)
       \ge S_0(t_0)e^{\tau G}S_0(t_0).
\]
The outer blocks retain no jumps but keep cooperative transfer.
In the middle block diagonal convolution commutes with transport;
the $N$ active jumps in \eqref{eq:active} are integrated over the
ordered simplex of volume $\tau^N/N!$, rather than prescribed at
single instants. At $T=2t_0+\tau$, Lemma~\ref{le:scan} gives
\begin{equation}\label{eq:Doeblin}
 S(T)f\ge\beta\one\sum_j\int_0^Lf_j,
 \qquad \beta=k_0^2\gamma_\tau>0.
\end{equation}
The first factor on the input side transfers total mass into species $a$;
the middle term makes it spatially constant there; the last scan reaches
all terminal species. No characteristic is shared between different species.

\paragraph{Principal normalization.}
On the disjoint union of $m$ circles define
\begin{equation*}
 (P_tf)_i=e^{-st}\phi_i^{-1}(S(t)(\phi f))_i,\qquad
 \pi(f)=\sum_i\int_0^L\psi_i\phi_i f_i\,dx,\qquad
 \nu(f)=\frac1{mL}\sum_i\int_0^Lf_i.
\end{equation*}
The constructed principal pair makes $P_t$ Markov, meaning positive
and constant-preserving, and $\pi$ an invariant probability functional. From \eqref{eq:Doeblin},
\[
 P_Tf\ge\alpha\nu(f)\one\quad(f\ge0),\qquad
 \alpha=\min\left\{\frac12,e^{-sT}\beta
                          \frac{\phi_{\min}}{\phi_{\max}}mL\right\}>0.
\]
For real $f$, write
$\osc(f)=\max_{i,x}f_i(x)-\min_{i,x}f_i(x)$, with extrema over all
components and $x\in\T_L$. The operator
$(P_T-\alpha\one\nu)/(1-\alpha)$ is positive and preserves constants.
Thus $\osc(P_Tf)\le(1-\alpha)\osc(f)$.
The invariant mean lies between its minimum and maximum, giving
\[
 \norm{P_{nT}f-\pi(f)\one}\le2(1-\alpha)^n\norm f.
\]
Intermediate times are handled by Markov contraction; complexification
costs at most another factor two. Applying this to $f/\phi$ proves
Theorem~\ref{th:mix}, with for example
\begin{equation}\label{eq:mix-constants}
 \kappa=-T^{-1}\log(1-\alpha),\qquad
 C_{\mathrm{mix}}=\frac{4\phi_{\max}}{(1-\alpha)\phi_{\min}}.
\end{equation}

\subsection{Perturbations, phases and velocity derivatives}\label{sec:perturb}

Fix finitely many directed paths witnessing strong connectivity. If the entrywise
coefficient perturbation is smaller than $I_{ij}/(2L)$ on every used edge,
at least half of each required integral remains. Taking also
$\Delta_0<c_-/2$ preserves all nonzero velocities. The principal theorem
therefore applies on the same domain $D$ throughout this neighborhood.
For its fixed positive $\phi$,
\begin{equation*}
 \norm{(\widetilde A-A)\phi}\le C_0\Delta,\qquad
 |\widetilde s-s|\le C_0\Delta/\phi_{\min}.
\end{equation*}
The first estimate uses $(\widetilde C-C)\phi'$, the norm-one convolution
bound and $\norm{\widetilde K_i-K_i}\le\norm{\widetilde J_i-J_i}_1$;
the second tests both perturbed formulas \eqref{eq:CW} on $\phi$.

Choose one real $\lambda>0$ above all row bounds of $B,\widetilde B$
and their transposes. The Laplace formula gives common resolvent bounds.
The displayed resolvent equation gives common $X\to C^1$ bounds.
The common-domain resolvent identity is consequently a bounded composition:
\begin{equation*}
 \widetilde R_\lambda-R_\lambda
     =\widetilde R_\lambda(\widetilde A-A)R_\lambda,
 \qquad \norm{\widetilde R_\lambda-R_\lambda}\le C_1\Delta.
\end{equation*}
The derivative in $\widetilde A-A$ is applied only to $R_\lambda X\subset D$.
There is no norm-continuity assertion for the translation groups themselves.

Let $r=(\lambda-s)^{-1}$. Choose a small conjugation-symmetric circle
$\Gamma$ around this simple isolated resolvent eigenvalue, excluding zero
and every other spectral point. Uniformly on $\Gamma$, for small $\Delta$,
\[
 (z-\widetilde R_\lambda)^{-1}
 =(z-R_\lambda)^{-1}\sum_{n\ge0}
       [(\widetilde R_\lambda-R_\lambda)(z-R_\lambda)^{-1}]^n.
\]
The contour projections $E,\widetilde E$ differ by $O(\Delta)$ in norm.
Projections at distance less than one have isomorphic finite-dimensional
ranges, so $\widetilde E$ has rank one. The bound on $|\widetilde s-s|$ places $(\lambda-\widetilde s)^{-1}$ within $\Gamma$, identifying this
range as the perturbed principal eigenspace. Therefore
\begin{equation*}
 \widetilde\phi=\frac{\widetilde E\phi}{\mathcal N(\widetilde E\phi)},
 \qquad \norm{\widetilde\phi-\phi}\le C_2\Delta.
\end{equation*}
The denominator tends to one. The projection is real, and its normalized
vector is positive by closeness to $\phi\gg0$. In particular, it has a
common positive floor and upper bound.

The individual eigen-equations give
\begin{equation*}
 \phi_i'=c_i^{-1}\left[s\phi_i-d_i(K_i\phi_i-\phi_i)
                                      -\sum_jB_{ij}\phi_j\right].
\end{equation*}
Their tilded differences are $O(\Delta)$ in $X$: reciprocal velocities
are Lipschitz away from zero and all other terms are already controlled.
The normalized eigenvector estimate therefore improves to a $C^1$ bound. Apply the
same argument to the separate formal-adjoint realizations, first normalizing
their positive vectors $\chi,\widetilde\chi$ by $\mathcal N=1$.
Their positive pairings with $\phi,\widetilde\phi$ stay away from zero
and are Lipschitz. Division by the respective pairings gives the required
$C^1$ bound for $\psi,\widetilde\psi$. This proves \eqref{eq:parameter}.

The branch is in fact continuously differentiable in all the stated
data norms, including the $C^1$ norms of both principal functions.
To check this without a differentiability assumption on spatial
translations, regard $A$ as a bounded map $D\to X$ and fix the above
real resolvent point. The map from $(c,d,B,J)$ to this operator is
polynomial and hence analytic on the ambient Banach space. For small
unrestricted real perturbations, with $E=\widetilde A-A$,
\[
 (\lambda-\widetilde A)^{-1}
 =R_\lambda(\Id-E R_\lambda)^{-1}
\]
is analytic in $\mathcal L(X,D)$ by its norm-convergent Neumann series.
The same statement holds for the separate adjoint realization.
Contour calculus on $X$ makes the simple resolvent eigenvalue and its
rank-one projection analytic. If $w$ is its projected normalized
vector and $r$ the simple eigenvalue, the identity
$w=r^{-1}(\lambda-\widetilde A)^{-1}w$ upgrades the branch to an
analytic $D$-valued branch. The normalization denominators are nonzero
near the fixed pair. This proves the claimed local $C^1$ extension;
on the admissible set it is the positive branch already identified.
It also makes the derivative in \eqref{eq:HF} continuous.

Fix an active species $a$. In the proof of Lemma~\ref{le:convolution},
choose a bounded subdensity $q=a_01_E\le j$, where $j$ is the fixed
periodization of $J_a$. Suppose $q*q\ge b_0>0$ on an arc $I$.
Periodization is $L^1$ contractive on differences. With
$\widetilde q=\min(q,\widetilde j)$,
\begin{equation*}
 0\le\widetilde q\le a_0,\quad
 \norm{q-\widetilde q}_1\le\Delta,\quad
 \norm{q*q-\widetilde q*\widetilde q}_\infty\le2a_0\Delta.
\end{equation*}
The last bound follows by expanding the difference and applying the
$L^1$--$L^\infty$ convolution inequality. If
$\Delta<b_0/(4a_0)$, then $\widetilde j^{*2}\ge b_0/2$ on $I$.
A fixed number of arc convolutions therefore gives common $N$ and
$\epsilon_0>0$ with
\begin{equation*}
 \widetilde K_a^Nf\ge\epsilon_0\int_0^Lf\qquad(f\ge0).
\end{equation*}
The bounded function here is $q$, not either original kernel density.

Take also $\Delta_0<d_a/2$, so $\widetilde d_a\ge d_a/2$.
There are common bounds $0<c_*\le|\widetilde c_i|\le c^*$ and a common
diagonal compensation $H_0^*$. For $m>1$, let $I_*>0$ be the minimum
integral among the finitely many fixed used edges, so each perturbed used
edge has integral at least $I_*/2$. Set
\[
 \eta_* =\begin{cases}\min\{1,I_*/(2c^*)\},&m>1,\\1,&m=1,\end{cases}
 \quad
 t_* =mL/c_*,\quad k_*=e^{-H_0^*t_*}\eta_*^{m-1}.
\]
The proof of Lemma~\ref{le:scan} applies on
the retained strongly connected subgraph and its transpose: each destination cycle
has duration at most $L/c_*$ and each used edge factor is at least $\eta_*$.
Hence both scans have the common constant $k_*$ at $t_*$.
At $T_*=2t_*+1$, equations \eqref{eq:active}--\eqref{eq:Doeblin} give
\[
 \widetilde S(T_*)f\ge\beta_*\one\sum_j\int_0^Lf_j,
 \qquad
 \beta_*=k_*^2e^{-H_0^*}\frac{(d_a/2)^N}{N!}\epsilon_0>0.
\]
The already proved eigenfunction estimates give a common ratio bound
$\widetilde\phi_{\max}/\widetilde\phi_{\min}\le R_*$ and a common upper
bound $\widetilde s\le s_*$. Thus the transformed Markov semigroups have
the common minorization mass
$\alpha_* =\min\{1/2,e^{-s_*T_*}\beta_*mL/R_*\}>0$.
The oscillation argument gives common mixing constants, including all
intermediate times. The common time need not be a multiple of any cycle.

\subsubsection{Spatial phases and parameter derivatives}
Simultaneous translation $(T_\vartheta f)(x)=f(x+\vartheta)$ commutes
with every $c_i\partial_x$ and every $K_i$. Therefore
\[
 A_\vartheta=T_\vartheta A T_{-\vartheta},\qquad
 s_\vartheta=s,\qquad
 \phi_\vartheta=T_\vartheta\phi,\quad
 \psi_\vartheta=T_\vartheta\psi.
\]
All normalizations, contour norms, derivative bounds, graph integrals,
kernel minorants and eigenfunction ratios used above are unchanged by
this isometric conjugacy. This proves phase uniformity. Uniform-norm differentiability under
translation would require additional regularity of the coefficients.

For a differentiable admissible path, the established $C^1$
dependence into $D$ allows differentiation of $A\phi=s\phi$ in $X$.
Pairing with $\psi$ and using \eqref{eq:pairing} cancels the
$\dot\phi$ terms and gives \eqref{eq:HF}. In particular the velocity
contribution is $\sum_i\dot c_i\int\psi_i\phi_i'$.
This completes the proof of Theorem~\ref{th:perturb}.

The autonomous parameter map is defined on an ambient Banach space
with $c_i\ne0$, coefficient topology $C(\T_L)$ and kernel topology
$L^1(\R)$. The mass-one equations $\int J_i=1$ define affine
hyperplanes; admissible probabilities are their intersections with the
nonnegative cone. The latter constraint is not an open affine manifold.
Thus an admissible $C^1$ path satisfies $\int\dot J_i=0$, as well as
its positivity constraints, while the local analytic projection is
constructed in the ambient open operator neighborhood.
If $R(z)=(z-A)^{-1}$ and $\Gamma$ is the fixed small positively oriented
contour about the simple principal value, its derivative is
\begin{equation*}
 \dot\Pi=\frac1{2\pi i}\int_\Gamma R(z)\dot A R(z)\,dz,
 \qquad \dot A\in\mathcal L(C^1,X).
\end{equation*}
Indeed $R:X\to C^1$ is uniformly bounded on $\Gamma$; the fixed-domain
Neumann identity gives $\dot R=R\dot A R$ uniformly there, so the
finite contour integral can be differentiated in operator norm.
This also gives its $X\to C^1$ interpretation. The coefficient
derivative uses the uniform norm and the kernel derivative the $L^1$
norm. It includes the autonomous terms $\dot c_i\partial_x$.

For periodic coefficients the velocities remain fixed. The domain is
the characteristic domain $\mathcal D$ of
Proposition~\ref{tp:th:general}, with continuous derivatives along
transport directions. Unlike the autonomous $C^1$ adjoint realization,
periodic duals may be Borel measures. Translation is not operator-norm
continuous on $C(\T_L)$, so an operator-norm Duhamel argument does not
yield velocity differentiation for these rough periodic coefficients.

\section{Periodic Floquet theory and characteristic communication}\label{sec:periodic-theory}

The periodic setting calls for a different argument: velocities may
vanish, so the autonomous compact-resolvent theory does not apply. We
compare the full period map with its no-jump part; compactness of their
difference underlies the resulting spectral criterion.

\subsection{Compact splitting and principal normalization}\label{sec:compact-return}
We collect the standard positive-operator and compact-perturbation
steps with a quantitative Doeblin estimate. For the transport application,
$E$ is the disjoint union of the component circles and
$\mu_0(f)=\sum_i\int_0^L f_i(x)\,dx$. Full support makes this mass
positive for every nonzero nonnegative continuous input.

\begin{theorem}[Positive compact splitting and principal normalization]
\label{thm:compact-return}
Let $E$ be a compact metric space, let $\mu_0$ be a finite nonzero
Borel measure with full support, and put $X=C(E;\mathbb R)$.
Here $\phi\gg0$ means $\min_E\phi>0$. Let $M,V\ge0$ be bounded and
$M=V+Q$ with $Q\ge0$ compact. Suppose that for some positive
integers $N,R$ and $\eta,\delta>0$,
\[
 M^Nf\ge\eta\one\mu_0(f)\quad(f\ge0),\qquad
 (M^R-V^R)\one\ge\delta\one.
\]
Then $M$ has a strictly positive eigenvector if and only if
$r(M)>r(V)$. If $r(M)>r(V)$, then $r(M)>0$ is an algebraically
simple eigenvalue and $r(M)^{-n}M^n$ converges in operator norm
exponentially to its positive rank-one principal projection.
Once $M\phi=r(M)\phi$ with $\phi\gg0$, the quantitative bound is
\begin{equation}\label{eq:return-gap}
 r(V)^R\le r(M)^R-\delta\frac{\min\phi}{\max\phi}.
\end{equation}
The equivalence is qualitative; \eqref{eq:return-gap} is a
profile-dependent quantitative bound. The row tests below give gap
margins directly from coefficients, and the minorant quantifies normalized
mixing. Strictly positive return is a separate hypothesis.
\end{theorem}
\begin{proof}
\emph{The positive spectral point.} The minorant gives $M^N\one\ge\eta\mu_0(E)\one$, hence $r:=r(M)>0$.
We first verify the standard positive-operator fact $r\in\sigma(M)$.
Suppose to the contrary that $r\in\rho(M)$. The positive resolvents
$(z-M)^{-1}$ for real $z>r$ then converge in operator norm to
$\mathcal R=(r-M)^{-1}\ge0$ as $z\downarrow r$.
Choose $0<\epsilon<r$ with $\epsilon\|\mathcal R\|<1$ and put
$z=r-\epsilon$. Only under this contrary assumption, factorization gives
\[
 (z-M)^{-1}
   =\mathcal R(\Id-\epsilon\mathcal R)^{-1}
   =\sum_{n=0}^\infty\epsilon^n\mathcal R^{n+1}\ge0.
\]
Thus $w=(z-M)^{-1}\one$ satisfies
$zw=Mw+\one\ge\one$ and
$Mw\le(z-1/\|w\|)w$. Since $w\ge z^{-1}\one$ is an order unit,
iteration in its equivalent norm gives
$r(M)\le z-1/\|w\|<r$, a contradiction. Thus $r\in\sigma(M)$.

\emph{A positive eigenvector implies separation.} If $M\phi=\rho\phi$ with $\phi\gg0$, order domination and testing
on $\phi$ give $\rho=r(M)$. Since $M\ge V$, the difference
$M^R-V^R$ is positive. Its action on $\phi$ is at least
$\delta\min\phi\,\one$. Thus
\[
 V^R\phi\le
 \left(r(M)^R-\delta\frac{\min\phi}{\max\phi}\right)\phi,
\]
and iteration proves strict separation. The coefficient is nonnegative
because the left side is nonnegative.

\emph{Separation produces a positive eigenvector.} Conversely suppose $r(M)>r(V)$. For $|z|>r(V)$,
\[
 z-M=(z-V)[\Id-(z-V)^{-1}Q].
\]
On the complexification of $X$, the connected open set
$\Omega=\{z:|z|>r(V)\}$ contains $r$. The family
$(z-V)^{-1}Q$ is norm analytic and compact, and the bracket is invertible
for large $|z|$. The analytic Fredholm alternative
\cite[Chapter~VII, \S1, Theorem~1.9]{Kato} makes the singular set
discrete in $\Omega$; hence $r$ is isolated in $\sigma(M)$.
Near each singular point, the spectral subspace of $(z-V)^{-1}Q$
near $1$ is finite dimensional and varies analytically. After an
analytic identification with a fixed subspace, the inverse there is
a finite matrix inverse; on the complementary subspace it is analytic.
This reduction gives a meromorphic inverse with finite-rank principal
Laurent coefficients.
Since $r\in\sigma(M)$, the resulting singularity of $(z-M)^{-1}$
at $r$ cannot be removable: a norm limit there would be a two-sided
inverse of $r-M$. Thus, for an integer $k\ge1$,
\[
 (z-M)^{-1}=\sum_{j=1}^k(z-r)^{-j}P_{-j}+\mathcal H(z),
 \qquad P_{-k}\ne0,
\]
where $\mathcal H$ is analytic near $r$ and each $P_{-j}$ has finite
rank. Comparing the leading Laurent coefficient in
$(z-M)(z-M)^{-1}=\Id$ gives $(M-r)P_{-k}=0$.
For real $z>r$, the Neumann resolvent is positive, so
$P_{-k}=\lim_{z\downarrow r}(z-r)^k(z-M)^{-1}\ge0$.
If $P_{-k}\one=0$, then
$|P_{-k}f|\le\|f\|P_{-k}\one=0$ for every real $f$, a contradiction.
Hence $\phi=P_{-k}\one$ is a nonzero nonnegative eigenvector.
The integral minorant gives
$r^N\phi=M^N\phi\ge\eta\one\mu_0(\phi)\gg0$,
so $\phi$ is an order unit.

\emph{Algebraic simplicity.} Positivity and $M\phi=r\phi$ imply
$\|M^n\|\le(\max\phi/\min\phi)r^n$.
If $(M-r)v=u\ne0$ and $(M-r)u=0$, then
$r^{-n}M^nv=v+nr^{-1}u$, contradicting this bound.
For a small positively oriented circle $\Gamma$ enclosing only $r$,
its Riesz projection is
\[
 \mathcal P_r=\frac1{2\pi i}\int_\Gamma(z-M)^{-1}\,dz=P_{-1}.
\]
Its range is the finite-dimensional generalized eigenspace
\cite[Chapter~III, \S6.5]{Kato}. Every nontrivial Jordan chain contains
the pair just excluded, so this range equals the eigenspace.
For a real eigenvector $g$, put $c=\max_{x\in E}g(x)/\phi(x)$.
The nonnegative eigenvector $c\phi-g$ has a zero; unless it vanishes
identically, its minorant is strictly positive and contradicts that zero.
Real and imaginary parts then give one-dimensionality over the complex
field, and hence algebraic simplicity.

\emph{Operator-norm mixing.} Finally $Tf=r^{-1}\phi^{-1}M(\phi f)$ is Markov. Put
\[
 \Phi=\mu_0(\phi),\qquad
 \nu(f)=\Phi^{-1}\mu_0(\phi f),\qquad
 \epsilon=\min\left\{\frac12,
          \frac{\eta\Phi}{r^N\max\phi}\right\}>0.
\]
The raw minorant gives $T^Nf\ge\epsilon\nu(f)\one$ for $f\ge0$.
Since $T^N\one=\one$, subtracting this term leaves a positive
operator preserving $(1-\epsilon)\one$. Its division by
$1-\epsilon$ is Markov, so oscillation contracts by $1-\epsilon$
every $N$ steps. Minima and maxima of $T^nf$ are nested and their
common limit is a positive linear functional $\widehat\ell(f)$.
For real $f$ the error is at most
$2\|f\|(1-\epsilon)^{\lfloor n/N\rfloor}$; splitting into real
and imaginary parts gives the complex operator-norm bound
\[
 \|r^{-n}M^n-\Pi\|
 \le4\frac{\max\phi}{\min\phi}
          (1-\epsilon)^{\lfloor n/N\rfloor},\qquad
 \Pi f=\phi\widehat\ell(\phi^{-1}f).
\]
The limit is $\mathcal P_r$: it fixes $\phi$ and vanishes on
$\ker\mathcal P_r$, where $r-M$ is invertible and $\Pi(r-M)=0$.
The estimate specifies the profile bounds and minorant constants
that must remain uniform in a family.
\end{proof}

\subsection{A compact-perturbation Floquet criterion for transport}\label{sec:periodic}

For time-periodic coefficients, velocities may vanish. We first assume
that all jump rates are positive.

Fix spatial and temporal periods $L,T>0$ and $m\ge1$. On
$X=C(\T_L;\R^m)$ use the component maximum norm. Let
$c_i\in\R$ be fixed constants, with no nonzero or sign restriction,
and let $d_i>0$. Let $J_i$ be arbitrary nonnegative $L^1(\R)$
probability densities and
\[
 K_i f(x)=\int_\R J_i(h)f(x-h)\,dh,\qquad K=\operatorname{diag}(K_i).
\]
Let $B(t,x)$ be a continuous real matrix, periodic in $t,x$ with the
stated periods, and Metzler: $B_{ij}\ge0$ for $i\ne j$. Give the species
the edge $j\to i$ when
\[
 I_{ij}:=\int_0^T\int_0^L B_{ij}(t,x)\,dx\,dt>0.
\]
Assume this union graph is strongly connected, vacuously for $m=1$.
Set
\begin{equation*}
 \alpha=\min_{t,x,i}\sum_jB_{ij}(t,x),\qquad
 \beta=\max_{t,x,i}\left(\sum_jB_{ij}(t,x)-d_i\right).
\end{equation*}

Consider the mild periodic evolution
\begin{equation}\label{tp:eq:evolution}
 u_t=C\partial_xu+D_d(K-\Id)u+B(t,x)u,\qquad
 C=\operatorname{diag}(c_i),\quad D_d=\operatorname{diag}(d_i).
\end{equation}
Write $U(t,s)$ for its evolution family and $M=U(T,0)$.
Let $V(t,s)$ be the evolution with the positive jump terms $D_dK$
deleted, retaining $C\partial_x+B-D_d$, and put $V_0=V(T,0)$.
The spectral separation below is the usual positive compact-perturbation
criterion; compare \cite[Theorem~2.13 and Remark~2.14]{FLRX},
\cite[Theorem~2.1]{BS}, \cite[Lemmas~2.1--2.3]{WZ}, and
\cite{WLSV} in their no-transport settings.
The proof here verifies compactness and the chronological positivity
estimate for arbitrary $L^1$ circle kernels and component velocities.
The generator domain uses differentiation along component characteristics
(Appendix~\ref{app:periodic-framework}). The spatial duals are finite Borel
measures on $\T_L$, without an assumed continuous density.

\begin{theorem}\label{tp:th:main}
Under the preceding coefficient, kernel and graph assumptions, the
following conditions are equivalent:
\begin{enumerate}
\item There are a real $\lambda$ and a strictly positive continuous
      space-time periodic mild Floquet function $p$, so
      $U(t,s)p(s)=e^{\lambda(t-s)}p(t)$.
\item The full and no-jump period maps satisfy
\begin{equation}\label{tp:eq:spectral-gap}
 r(M)>r(V_0).
\end{equation}
\end{enumerate}
When these conditions hold, the spectral radius $r(M)>0$ is an
algebraically simple eigenvalue with a strictly positive eigenvector $\phi$.
With $\lambda=T^{-1}\log r(M)$ there is a strictly positive continuous
space-time periodic function
$p(t)=e^{-\lambda t}U(t,0)\phi$. It is unique up to a positive scalar.
There are positive bounded linear functionals $\ell_s$ on $X$, periodic
in $s$, with
\[
 \ell_s(p(s))=1,\qquad
 \ell_tU(t,s)=e^{\lambda(t-s)}\ell_s,
\]
and constants $A,\kappa>0$ such that
\begin{equation}\label{tp:eq:mix}
 \|e^{-\lambda(t-s)}U(t,s)f-p(t)\ell_s(f)\|
 \le A e^{-\kappa(t-s)}\|f\|\qquad(t\ge s).
\end{equation}
The actual uniform evolution growth bound equals $\lambda$, and
$\lambda\ge\alpha$.

The perturbation and first-variation results of
Section~\ref{sec:periodic-perturb} fix $L,T,c_i$: $B$ varies in the
uniform norm, the rates in Euclidean norm and kernels in $L^1$.
Admissible kernels remain nonnegative with integral one, reactions remain
Metzler and required active rates remain positive.
Profile derivatives are understood in the characteristic direction;
the dual functionals are represented by finite positive Borel measures.
\end{theorem}

\begin{corollary}\label{tp:cor:row}
If $\alpha>\beta$, then \eqref{tp:eq:spectral-gap} holds and all
conclusions of Theorem~\ref{tp:th:main} follow, with $\lambda>\beta$.
In particular, this applies for every positive set of jump rates when
all rows of $B(t,x)$ have the same constant sum.
\end{corollary}

All operator integrals below are vectorwise strong Bochner integrals;
the uniform estimates give the operator-norm bounds and parameter
derivatives.

\subsubsection{Generation and compact perturbation}
The diagonal translations
$(\mathcal T(t)f)_i(x)=f_i(x+c_i t)$ form a strongly continuous isometric
group. All remaining operators in \eqref{tp:eq:evolution} are bounded,
with norm-continuous time dependence. Iterating the transport Duhamel
formula gives a unique strongly continuous evolution on each finite
time interval: the $n$th term is bounded by
$C^n(t-s)^n/n!$. The series converges uniformly on each such interval
on bounded input sets. Uniqueness follows by Gronwall, and the
evolution property follows by uniqueness.

Choose $H>0$ such that
$B_{ii}(t,x)-d_i+H\ge0$ for every $i,t,x$.
The base group $e^{-H(t-s)}\mathcal T(t-s)$ and the positive
perturbations $D_dK$, the off-diagonal entries of $B$, and
$\operatorname{diag}(B_{ii}-d_i+H)$ give a termwise positive Dyson
series. In particular $U$ is positive. The same reasoning applies to
the no-jump evolution $V(t,s)$ generated by
$C\partial_x+B(t,x)-D_d$. Comparison with constant vectors gives
\begin{equation}\label{tp:eq:bounds}
 U(t,s)\one\ge e^{\alpha(t-s)}\one,\qquad
 \|U(t,s)\|\le e^{b_+(t-s)},\qquad
 \|V(t,s)\|\le e^{\beta(t-s)},
 \quad b_+=\max_{t,x,i}\sum_jB_{ij}(t,x).
\end{equation}
These inequalities follow first by inserting the constant
sub/supersolutions in the integral equation. Positivity then gives
the norm bounds. Periodicity implies
$M^n=U(nT,0)$ and $V(T,0)^n=V(nT,0)$.

Each $K_i$ is compact on $C(\T_L)$. Indeed periodize $J_i$ to
$j_i\in L^1(\T_L)$; Tonelli gives $\int_0^Lj_i=1$ and
\[
 \sup_{\|f\|\le1}\|K_i f(\cdot+h)-K_i f\|_\infty
 \le \|j_i(\cdot+h)-j_i\|_1\longrightarrow0.
 \]
Translation continuity in $L^1$, boundedness and Arzel\`a--Ascoli
prove the assertion. Therefore $D_dK$ is compact on $X$.

\begin{lemma}\label{tp:lem:compact}
For every fixed $t>s$, $U(t,s)-V(t,s)$ is compact.
\end{lemma}
\begin{proof}
The bounded-perturbation Dyson identities, or direct substitution and
Fubini in their absolutely bounded series, give pointwise Bochner
integrals
\[
 U(t,s)-V(t,s)=\int_s^t U(t,r)D_dK V(r,s)\,dr .
 \]
Let $\|f\|\le1$. The set of $D_dKV(r,s)f$ over all $r\in[s,t]$
is contained in the image under the compact operator $D_dK$ of one
fixed bounded ball. Its closure is compact. Joint continuity of
$(r,z)\mapsto U(t,r)z$ on that compact set times $[s,t]$ follows
from strong continuity and a common finite-interval operator bound.
Thus all integrand values are contained in one compact subset
$\mathcal C$ of $X$. Every integral divided by $t-s$ belongs to the
closed convex hull of $\mathcal C$, which is compact in a Banach space.
A finite $\varepsilon$-net of $\mathcal C$ approximates
its convex hull within $\varepsilon$ by the compact convex hull of
finitely many points; completeness gives compactness of the closure.
This proves collective relative compactness of the image of the unit
ball. No operator-norm continuity of the translation group or
interchange with an operator-norm integral was assumed.
\end{proof}

It follows that
\begin{equation*}
 M=V_0+Q,\quad Q\ \hbox{compact},\qquad
 r(V_0)\le e^{\beta T},\qquad e^{\alpha T}\le r(M).
\end{equation*}
The last inequality follows by iterating the first comparison in
\eqref{tp:eq:bounds}. This proves Corollary~\ref{tp:cor:row} once
the criterion is established. No row gap is used for compactness.

\subsubsection{Generalized periodic growth without a Floquet function}
The realization and growth formulas below remain valid when
\eqref{tp:eq:spectral-gap} fails. The spectrum correspondence follows
from \cite[Propositions~3.1--3.2]{HR}; the strict-stability equivalences
are instances of \cite[Theorems~3.1 and~4.7]{GM}. Appendix~\ref{app:periodic-framework}
gives a direct verification, including the characteristic domain and the
surjectivity needed for the reverse resolvent implication. Write
$\T_T=\R/T\mathbb Z$ for the time circle. On
$Y=C(\T_T\times\T_L;\R^m)$ let $\mathcal D$ be the generator domain of
the strongly continuous group
\[
 (\mathcal H(\tau)f)_i(t,x)=f_i(t-\tau,x+c_i\tau).
\]
Equivalently, each component has a continuous derivative in the indicated
characteristic direction. The fundamental theorem along a characteristic
and compactness give uniform difference quotients. Separate
$\partial_t$ and $\partial_x$ are not required. Put
\[
 \mathcal B(t)=D_d(K-\Id)+B(t,x),\qquad
 \mathscr G=(C\partial_x-\partial_t)+\mathcal B(t),\qquad
 D(\mathscr G)=\mathcal D,
 \qquad \lambda_0=T^{-1}\log r(M).
\]
Here $r(M)>0$ follows from \eqref{tp:eq:bounds}.
The directional derivative has the group meaning just defined.
Bounded perturbation gives the positive evolution semigroup
\begin{equation*}
 (\mathcal E(\tau)f)(t)=U(t,t-\tau)f(t-\tau),\qquad \tau\ge0,
\end{equation*}
as is seen by matching the two Dyson series.

The uniform evolution growth bound $\omega_{\mathrm{unif}}(U)$ is the
infimum of $\omega$ for which
$\|U(t,s)\|\le C_\omega e^{\omega(t-s)}$ for every $t\ge s$.

\begin{proposition}\label{tp:th:general}
On the complexifications,
\begin{equation}\label{tp:eq:spectrum-general}
 z\in\rho(\mathscr G)\quad\Longleftrightarrow\quad
 e^{zT}\in\rho(M).
\end{equation}
Consequently, without a positive Floquet function,
\begin{equation}\label{tp:eq:growth-general}
 s(\mathscr G)=\omega_0(\mathcal E)
       =\omega_{\mathrm{unif}}(U)=\lambda_0.
\end{equation}
Moreover,
\begin{equation}\label{tp:eq:upper-general}
 \lambda_0=\inf_{\substack{q\in\mathcal D\\\min q>0}}
                  \max_{t,x,i}\frac{(\mathscr Gq)_i}{q_i}.
\end{equation}
The following are equivalent:
\begin{equation*}
 \lambda_0<0;\qquad
 \exists q\in\mathcal D:\ \min q>0,\ \mathscr Gq\ll0;\qquad
 0\in\rho(\mathscr G),\quad(-\mathscr G)^{-1}\ge0.
\end{equation*}
Strict stability implies the periodic maximum principle on $\mathcal D$.
Without a positive Floquet function, no converse from that plain maximum
principle at the zero threshold is asserted.
\end{proposition}
A direct proof is supplied in Appendix~\ref{app:periodic-framework}.

\subsubsection{Communication estimates and Floquet normalization}
Apply Lemma~\ref{le:convolution} to each component kernel. Its proof
works with any bounded nonzero subdensity, in particular
$q_i=\min(j_i,A_i)$ for a fixed $A_i>0$.

Let $n_i,\eta_i$ be such that $K_i^{n_i}f\ge\eta_i\int f$
for $f\ge0$. Retain only $n_i$ jumps in component $i$ in its
diagonal Dyson series over a period. Translations commute with
that component's convolution, so for every starting time $s$
\begin{equation}\label{tp:eq:onejump}
 (U(s+T,s)f)_i\ge
 \zeta_i\int_0^L f_i,\qquad
 \zeta_i=e^{-HT}(d_iT)^{n_i}\eta_i/n_i!>0 .
\end{equation}
This lower bound ignores all switching and all other positive
Dyson terms. It remains true for zero or negative $c_i$.

Suppose a component $j$ starts with a spatial constant floor $b>0$.
Retaining one off-diagonal multiplication and no other positive
terms over the next period gives
\[
 (U(s+T,s)f)_i(x)\ge
 b e^{-HT}\int_s^{s+T}B_{ij}(r,x+c_i(s+T-r))\,dr.
\]
Integration in $x$ and temporal periodicity give target mass at least
$b e^{-HT}I_{ij}$. A further period and \eqref{tp:eq:onejump} turn
this mass into a spatial constant floor at least
$b e^{-HT}I_{ij}\zeta_i$.

For each ordered pair choose a directed path of length
$l\le m-1$, taking length zero for equal indices. One initial jump
period and $l$ transfer/jump pairs give a positive floor from its
input mass after $(1+2l)T$. Pad the remaining periods using the
zero-jump lower bound $e^{-H\Delta}\mathcal T(\Delta)$ on constants.
With $N=2m-1$, the finite minimum of all path constants gives
\begin{equation}\label{tp:eq:allminorant}
 (U(s+NT,s)f)_i(x)\ge
 \zeta_*\sum_j\int_0^L f_j(y)\,dy
 \quad(f\ge0,\ s\in\R),\qquad \zeta_*>0 .
\end{equation}
For $m=1$, $N=1$ and only \eqref{tp:eq:onejump} is used.
All paths are chronological; no simultaneous pointwise
irreducibility of $B$ is required.

Let $d_-:=\min_i d_i>0$. Apply zero-jump retention directly to
the fixed test vector $\one$:
$V(t,0)\one\ge e^{-Ht}\one$ and
$U(T,t)\one\ge e^{-H(T-t)}\one$. Since $K_i\one=\one$,
\begin{equation*}
 (M-V_0)\one
 =\int_0^T U(T,t)D_dK V(t,0)\one\,dt
 \ge T d_-e^{-HT}\one.
\end{equation*}
This return estimate precedes the construction of any principal vector.

Theorem~\ref{thm:compact-return}, applied to the compact splitting with
input $\one$ and using \eqref{tp:eq:allminorant} together with the strict
lower bound above, proves both directions of the gap criterion, algebraic
simplicity and discrete mixing. Normalize its vector $\phi$ by
$(mL)^{-1}\sum_i\int\phi_i=1$, set $r=r(M)$, $\lambda=T^{-1}\log r$,
and $p(t)=e^{-\lambda t}U(t,0)\phi$. The eigenvector identity and
periodicity give a continuous periodic mild order unit. Define
$p_-:=\min_{t,x,i}p_i(t,x)>0$ and
$p_+:=\max_{t,x,i}p_i(t,x)<\infty$, with extrema over one period cell.
Its covariance implies
$M_sp(s)=rp(s)$, $M_s=U(s+T,s)$, and order domination gives $r(M_s)=r$.

In \eqref{tp:eq:allminorant}, $N=2m-1$ and $\zeta_*>0$ are
independent of the starting phase. Set
$\epsilon:=\min\{1/2,\zeta_*r^{-N}mLp_-/p_+\}>0$.
The Markov maps $Q_sf=r^{-1}p(s)^{-1}M_s(p(s)f)$ obey
\[
 Q_s^Nf\ge\epsilon\nu_s(f)\one,\qquad
 \nu_s(f)=\frac{\sum_i\int p_i(s,x)f_i(x)\,dx}
                   {\sum_i\int p_i(s,x)\,dx}.
\]
The subtraction, oscillation contraction and nested extrema proved in
Theorem~\ref{thm:compact-return} apply with this common $\epsilon$.
Write their limit as $\gamma_s(f)$ and put
$\ell_s(f)=\gamma_s(f/p(s))$. This gives the unique positive left
eigenfunctional normalized by $\ell_s(p(s))=1$. Evolution and
periodicity imply $U(T,t)M_t=M_0U(T,t)$, so uniqueness yields
\begin{equation}\label{tp:eq:dual}
 \ell_t=e^{-\lambda(T-t)}\ell_0U(T,t)\qquad(0\le t\le T).
\end{equation}
Strong continuity gives weak-star continuity and the asserted covariance.
Periodic extension and a residual time interval give \eqref{tp:eq:mix}:
the residual evolution is Markov in the $p$ weights, so cannot increase
oscillation and has original norm at most $p_+/p_-$. Testing on $p(s)$
and order domination give the two-sided growth bounds and identify
$\lambda$ as the actual uniform growth bound. A second positive Floquet
function gives a positive eigenvector of $M$, hence a scalar multiple.

In the preceding space-time realization, the identity $\mathcal E(\tau)p=e^{\lambda\tau}p$ implies
$p\in\mathcal D$ and $\mathscr Gp=\lambda p$ by taking the norm
derivative at zero. Its order-unit bounds show that both the
spectral and growth bounds of $\mathscr G$ equal $\lambda$.

For $\mathcal P=\{w\in\mathcal D:\min w>0\}$, positive comparison
with $e^{\mu\tau}w$ yields
\begin{equation}\label{tp:eq:cw}
 \lambda=
 \sup_{w\in\mathcal P}\min_{t,x,i}
       \frac{(\mathscr Gw)_i}{w_i}
 =
 \inf_{w\in\mathcal P}\max_{t,x,i}
       \frac{(\mathscr Gw)_i}{w_i}.
\end{equation}
For clarity, $\mathscr Gw\ge\mu w$ implies
$\mathcal E(\tau)w\ge e^{\mu\tau}w$ by its integral identity and
positivity, and hence $\mu\le\lambda$ by the actual growth bound.
The upper inequality is analogous. Both extrema in
\eqref{tp:eq:cw} are attained at $p$.

The following are equivalent:
\[
 \lambda<0;\qquad
 [f\in\mathcal D,\ \mathscr Gf\ge0\Longrightarrow f\le0];
 \qquad 0\in\rho(\mathscr G),\ (-\mathscr G)^{-1}\ge0 .
\]
The implication from strict stability follows from
Proposition~\ref{tp:th:general}. If $\lambda\ge0$, the positive domain
function $p$ violates the maximum principle; at $\lambda=0$ the inverse
fails, and at $\lambda>0$ a positive inverse would send $p$ to
$-p/\lambda$. This proves the converse implications.
These equivalences hold on the directional generator domain.
The space-time semigroup retains temporal phase harmonics, whereas
\eqref{tp:eq:mix} concerns the evolution between fixed phases.

\subsubsection{Fixed-velocity perturbations and phases}\label{sec:periodic-perturb}
For both periodic classes, fix
$L,T$ and velocities. Use
\[
 \Delta=\|\widetilde B-B\|_\infty+
 \max_i|\widetilde d_i-d_i|+
 \max_i\|\widetilde J_i-J_i\|_1.
\]
Kernels remain probabilities, reactions Metzler, and the active rates
positive. The finitely chosen graph integrals or tree scan margins
persist: their changes are bounded by $TL\Delta$ or $W\Delta$.
The transport group is unchanged, so bounded-perturbation Duhamel gives
$\|\widetilde U(t,s)-U(t,s)\|+\|\widetilde V(t,s)-V(t,s)\|\le C\Delta$
on each fixed period block.

Let $R$ be the return length in periods: $R=1$ here and $R=N_R$
in Section~\ref{sec:one-active}. Choose the raw return constant
$\delta>0$ with $(M^R-V_0^R)\one\ge\delta\one$.
For the fixed principal vector $\phi$, positivity gives
\[
 V_0^R\phi\le b_0^R\phi,\qquad
 b_0^R=r^R-\delta\phi_-/\phi_+<r^R,
 \quad \phi_-:=\min\phi,\quad\phi_+:=\max\phi.
\]
Here $b_0>0$ by no-jump retention. Perturbed order-unit bounds give
$r(\widetilde M)\ge r-C\Delta$ and
$r(\widetilde V_0)\le(b_0^R+C\Delta)^{1/R}\le b_0+C'\Delta$.
Thus the gap persists; no general Lipschitz assertion for spectral radii
is used. Simple contour projections give local Lipschitz principal data
$r,\lambda,\phi,\ell_0$. The mild profile and \eqref{tp:eq:dual} give
the same bounds for $p(t)$ and $\ell_t$ uniformly over a period,
in $C$ and dual norms. The directional equation gives graph-norm
dependence on $\mathcal D$ because its directional part is fixed.

The kernel constants are locally uniform. For each used active kernel,
take the bounded subdensity $q_i=\min(j_i,A_i)$ from
Lemma~\ref{le:convolution} and
$\widetilde q_i=\min(q_i,\widetilde j_i)$. Then
\[
 \|q_i-\widetilde q_i\|_1\le\Delta,\qquad
 \|q_i*q_i-\widetilde q_i*\widetilde q_i\|_\infty\le2A_i\Delta.
\]
The same arc, half its lower height and the same convolution power work.
Together with the persistent finite paths and positive profile bounds,
this supplies common communication and mixing constants.

For a differentiable admissible bounded-data path, writing
$\mathcal B(t)=D_d(K-\Id)+B(t,\cdot)$ gives
\[
 M'=\int_0^T U(T,t)\dot{\mathcal B}(t)U(t,0)\,dt,\qquad
 \dot\lambda=\frac1T\int_0^T\ell_t(\dot{\mathcal B}(t)p(t))\,dt.
\]
The first integral is strong. Dividing Duhamel by the parameter increment
$h$, its difference from that integral has operator norm
$CT(o(1)+O(|h|))$: the bounded-data quotient converges uniformly and
$\|\widetilde U-U\|=O(|h|)$. This proves the operator-norm derivative
without requiring an operator-norm integrand. Differentiating the simple
eigenvalue and using \eqref{tp:eq:dual} gives the second formula.

Time and space translations preserve all gate minima, graph integrals,
diagonal bounds and kernel constants. Their profile is the translated
$p$, up to normalization, with the same exponent and ratio $p_-/p_+$.
The return estimate therefore supplies a common gap margin, so the
neighborhoods and mixing constants above are uniform in both phases.
No phase derivative is claimed. Gate positivity also persists under small
velocity changes by uniform continuity of $B$, but this does not give an
operator-norm velocity perturbation theorem for rough periodic data.

\begin{remark}\label{rem:periodic-row-scope}
The row separation $\alpha>\beta$ is a sufficient dispersal-versus-
heterogeneity condition, not an assumption that a principal
eigenfunction already exists. For constant row sum $b$ it holds
for every $d_i>0$. In general it holds whenever all $d_i$ exceed
the oscillation between the smallest and largest row sums.
The resonant continuous-cusp example with
$B(t,x)=(1-a(x+t))\Id+
\left(\begin{smallmatrix}-1&1\\1&-1\end{smallmatrix}\right)$ and
$d_i=1$ has $\alpha=-15$, $\beta=0$, so it violates this precise
condition. No contradiction with that example or with the
autonomous nonzero-drift theorem is implied.
\end{remark}

\subsection{A characteristic test with a designated active periodic component}
\label{sec:one-active}
In \eqref{tp:eq:evolution}, allow $d_i\ge0$ and designate an index
$a$ with $d_a>0$. The velocities remain arbitrary
constants and the kernels remain arbitrary $L^1$ probability densities.
The following coefficient condition replaces the union-graph assumption
of Section~\ref{sec:periodic}. It distinguishes receiving a pointwise
floor from sending spatial mass. Geometric control along characteristics is used in kinetic theory
\cite{DHHM}, and positive spectral and contraction methods are developed
in \cite{CG,CEAG}. For the present system, we construct the two scans,
combine them in chronological order and verify the spectral criterion.

Choose an integer $n_0\ge1$ and put $W=n_0T$. For $i\ne j$ define
\begin{align}
 P_{ij}&=\min_{s\in[0,T],\ y\in\T_L}
       \int_0^W B_{ij}(s+r,y-c_i r)\,dr,\label{oa:P}\\
 Q_{ij}&=\min_{s\in[0,T],\ y\in\T_L}
       \int_0^W B_{ij}(s+r,y-c_j r)\,dr.\label{oa:Q}
\end{align}
Continuity and compactness give these minima. The \emph{two-tree
condition} means that each species can be reached from $a$ by a path
of edges with $P_{ij}>0$, and can reach $a$ by a path of edges with
$Q_{ij}>0$. Choose such paths without repeated vertices. Let $p_*>0$
and $q_*>0$ be the respective minima, over species, of the products
of the edge numbers on these chosen paths. Empty products are one;
for $m=1$ both constants are one and no tree condition is needed.
Every chosen path has length at most $m-1$.

Fix $h\ge0$ with $h+B_{ii}-d_i\ge0$ everywhere. Retaining
$e^{-h(t-s)}\mathcal T(t-s)$ as the base evolution makes all the
remaining Dyson terms nonnegative, for both $U$ and $V$. In
particular the no-jump retention lower bound holds even when some
velocities or jump rates vanish.

\begin{lemma}\label{oa:scans}
For $f\ge0$ and $i\ne j$,
\begin{equation*}
 (V(s+W,s)f)_i(x)\ge e^{-hW}\int_0^W
 B_{ij}(s+r,x+c_i(W-r))
 f_j(x+c_i(W-r)+c_jr)\,dr.
\end{equation*}
A floor $f_j\ge b$ therefore gives a floor $be^{-hW}P_{ij}$
in component $i$, whereas
\[
 \int_0^L(V(s+W,s)f)_i\ge e^{-hW}Q_{ij}\int_0^L f_j.
\]
Put $H_S=(m-1)W$ and let $e_a$ be the constant unit vector of
species $a$. Then
\begin{equation}\label{oa:access}
 V(s+H_S,s)e_a\ge e^{-hH_S}p_*\one,\qquad
 \int_0^L(V(s+H_S,s)f)_a
       \ge e^{-hH_S}q_*\sum_j\int_0^L f_j.
\end{equation}
\end{lemma}
\begin{proof}
Keep the single off-diagonal multiplication $B_{ij}$ in the
no-jump Dyson series. The translations before and after it give the displayed estimate. In the floor estimate set $y=x+c_iW$ in
the coefficient integral. For the mass estimate first integrate
in $x$ and set $z=x+c_i(W-r)+c_jr$. The coefficient becomes
$B_{ij}(s+r,z-c_jr)$. Tonelli and \eqref{oa:Q} give the stated
bound. Thus the receiving velocity in $P$ and the sending velocity
in $Q$ arise from different, explicit changes of variables.

Iterate one edge per block of length $W$ along each chosen path.
Pad a shorter path at its end by zero-jump retention in its terminal
component; translation preserves constants and Lebesgue integrals.
Every starting phase is allowed in the minima. For the mass inequality
apply the inward path separately to $f_je_j$ and then add over $j$
by linearity. This does not add different lower bounds for the same
input. These steps prove \eqref{oa:access}; for $m=1$ it is the
identity at duration zero.
\end{proof}

\begin{lemma}
\label{oa:minorant-return}
Choose $n\ge1$ and $\eta_a>0$ from Lemma~\ref{le:convolution} so that
$K_a^nf\ge\eta_a\int_0^L f$ for $f\ge0$. Put
\[
 H_D=(2m-1)W=N_DT,\quad H_R=mW=N_RT,
\]
and
\begin{equation}\label{oa:constants}
 \zeta=e^{-hH_D}p_*q_*\eta_a\frac{(d_aW)^n}{n!}>0,
 \qquad \delta=e^{-hH_R}d_aWp_*>0.
\end{equation}
Then, uniformly in the starting phase,
\begin{align}
 U(s+H_D,s)f&\ge\zeta\one\sum_j\int_0^L f_j
                    &&(f\ge0),\label{oa:minorant}\\
 [U(s+H_R,s)-V(s+H_R,s)]\one&\ge\delta\one.
                                            \label{oa:return}
\end{align}
Only the outward tree is needed for \eqref{oa:return}.
\end{lemma}
\begin{proof}
Use consecutive intervals of lengths $H_S,W,H_S$. First use the
mass inequality in \eqref{oa:access}. In the middle retain exactly
$n$ jumps in species $a$, with no switches or further positive
multiplications. That component is bounded below by
\[
 e^{-hW}\frac{(d_aW)^n}{n!}\mathcal T_a(W)K_a^n f_a
 \ge e^{-hW}\frac{(d_aW)^n}{n!}\eta_a\int_0^L f_a.
\]
Here only the active component's own translation and convolution
are commuted. Finally use the floor inequality in
\eqref{oa:access}. Positivity and the evolution identity justify
each chronological replacement; multiplying constants proves
\eqref{oa:minorant}. The convolution lemma uses a bounded
nonzero subdensity, so it imposes no pointwise regularity or
positive lower bound on $J_a$ itself.

For \eqref{oa:return}, start instead from $e_a\le\one$ and retain
one active jump in the first interval of length $W$. Since
$K_a\one=\one$, the resulting active floor is $e^{-hW}d_aW$.
Apply the outward paths during the remaining $H_S$ units. Every
selected chronological word contains a spatial jump. The no-jump
series consists of precisely the words containing none, so their
complement is a positive series for $U-V$. The selected terms therefore give the return
bound for the difference, in addition to a lower bound for $U$.
\end{proof}

\begin{remark}
\label{rem:communication-logic}
Under the hypotheses of the theorem below, the two trees give
\eqref{oa:minorant} and \eqref{oa:return}, whereas a positive Floquet
function exists exactly when $r(M)>r(V_0)$. Corollary~\ref{oa:row}
verifies this separate gap from coefficients. It cannot be inferred
from the trees: Corollary~\ref{cor:one-active-no-gap} has both trees,
a minorant and positive return, but equal radii and no Floquet function.
Propositions~\ref{prop:smooth-lock}, \ref{prop:dual-lock} and
\ref{prop:near-resonance} exhibit the two scan failures and nonresonant
success.
\end{remark}

\begin{theorem}[Periodic spectral theory under the two-tree test]
\label{oa:criterion}
Under the two-tree condition, which requires only $d_a>0$ for the
designated component (the other rates may be positive or zero), the
conclusions of Theorem~\ref{tp:th:main} hold: a strictly positive
periodic mild Floquet function exists if and only if
$r(M)>r(V_0)$. When it exists the principal multiplier is algebraically
simple, \eqref{tp:eq:mix} holds, and the dual functionals are strictly
positive on every nonzero nonnegative input. The attained formulas
\eqref{tp:eq:cw} and the maximum-principle equivalences on $\mathcal D$
hold as well. The integral minorant \eqref{oa:minorant} holds
independently of the spectral separation.
\end{theorem}
\begin{proof}
Compactness of $U-V$ in Lemma~\ref{tp:lem:compact} uses only compactness
of $D_dK$, boundedness of $V$ and strong continuity on compact sets;
it remains valid when some $d_i=0$. The periodic directional realization
in Appendix~\ref{app:periodic-framework} likewise needs only bounded
perturbation of the translation group. Lemma~\ref{oa:minorant-return}
provides the two independent estimates required by
Theorem~\ref{thm:compact-return}, with return length $N_R$. The latter
therefore proves the exact gap criterion, simplicity and discrete mixing.

The phase construction above applies with $N_D$ and $\zeta$ from
Lemma~\ref{oa:minorant-return} in place of $N$ and $\zeta_*$.
For the same Markov maps $Q_s$, profile bounds $p_-,p_+$ and
probabilities $\nu_s$, choose
$\epsilon:=\min\{1/2,\zeta r^{-N_D}mLp_-/p_+\}>0$.
Then $Q_s^{N_D}f\ge\epsilon\nu_s(f)\one$ for $f\ge0$.
Oscillation contracts by $1-\epsilon$ every $N_D$ periods, and the
invariant functional is at least $\epsilon\nu_s$. Thus the dual is
strictly positive. Identity \eqref{tp:eq:dual} and normalized residual
intervals give continuous-time mixing. Finally
$\mathcal E(\tau)p=e^{\lambda\tau}p$ places $p$ in $\mathcal D$;
the order-unit comparisons proving \eqref{tp:eq:cw} and the maximum
principle apply unchanged. No separate spatial derivative or adjoint
density is used.
\end{proof}

\begin{corollary}[A coefficient test for the gap]\label{oa:row}
Define
\[
 b_+(t)=\max_{i,x}\sum_jB_{ij}(t,x),\quad
 b_-(t)=\min_{i,x}\sum_jB_{ij}(t,x),\quad
 \Delta_b=\int_0^T(b_+-b_-)\,dt.
\]
Choose $h_+\ge0$ with
$h_+\ge\max_{i,t,x}\{d_i-B_{ii}(t,x)+b_+(t)\}$, and put
\[
 \eta_R=\min\{1/2,e^{-h_+H_R}d_aWp_*\}>0.
\]
If $N_R\Delta_b<-\log(1-\eta_R)$, then
\begin{equation}\label{oa:rowgap}
 r(V_0)\le e^{\int_0^T b_+}(1-\eta_R)^{1/N_R}
       <e^{\int_0^T b_-}\le r(M).
\end{equation}
In particular, if every row sum equals a time-only function $b(t)$,
the gap and normalized mixing always hold under the two-tree test, and
\begin{equation}\label{oa:explicit-p}
 \lambda=\bar b:=T^{-1}\int_0^T b(t)dt,\qquad
 p(t)=\exp\!\left(\int_0^t(b(r)-\bar b)dr\right)\one.
\end{equation}
\end{corollary}
\begin{proof}
Put $\widehat U(t,s)=e^{-\int_s^t b_+}U(t,s)$ and define
$\widehat V$ in the same way. Their coefficient matrix is
$\widehat B=B-b_+(t)\Id$, with the same off-diagonal scans and $p_*$.
The choice of $h_+$ gives $h_++\widehat B_{ii}-d_i\ge0$.
Apply Lemma~\ref{oa:minorant-return} to these normalized evolutions:
\[
 [\widehat U(s+H_R,s)-\widehat V(s+H_R,s)]\one
 \ge e^{-h_+H_R}d_aWp_*\one\ge\eta_R\one.
\]
Thus $h_+$ already accounts for normalization. The row bounds give
$\widehat U(t,s)\one\le\one$, and hence
$\widehat V(s+H_R,s)\one\le(1-\eta_R)\one$.
Since $H_R=N_RT$, positivity and the spectral-radius formula give
the first bound in \eqref{oa:rowgap}; lower row comparison gives the last.
The strict inequality is exactly $N_R\Delta_b<-\log(1-\eta_R)$.
Common row sums give $\Delta_b=0$ and the explicit mild solution
\eqref{oa:explicit-p}, whose exponent follows by order-unit comparison.
\end{proof}

\paragraph{Perturbations and phases for partial dispersal.}
Section~\ref{sec:periodic-perturb} applies with the return length $N_R$
and the constants in \eqref{oa:constants}. The chosen gate integrals
change by at most $W\|\widetilde B-B\|_\infty$, so their positive
margins persist. The only required positivity is $d_a>0$ for the designated component; other rates may also be positive or may vanish.
Time and space phase shifts leave these scan minima unchanged. Thus
the fixed-velocity neighborhoods, first variation and mixing estimates
there are uniform in both phases for this class also.

\paragraph{Traveling gates restore communication away from resonance.}
Take $L=T=1$, two species, $d_1=d>0$, $d_2=0$, $\alpha_0>0$,
and a nonnegative nonzero smooth function $g$ on $\T$. Let
$\omega\in\mathbb Z$ and let $b(t)$ be continuous and one-periodic.
Set
\begin{equation*}
 B(t,x)=b(t)\Id+
 \begin{pmatrix}-\alpha_0&\alpha_0\\
       g(x+\omega t)&-g(x+\omega t)\end{pmatrix}.
\end{equation*}
The velocities $c_1,c_2$ may vanish or have opposite signs.
If $c_2\ne\omega$, choose an integer $W$ with
$|\omega-c_2|W\ge1$. Writing $G_0=\int_0^1 g>0$, direct integration
along the inactive characteristic gives
\[
 P_{21}\ge G_0/|\omega-c_2|>0,\qquad Q_{12}=\alpha_0W>0.
\]
At least one full circuit occurs, for either sign of the relative
velocity. These are the two required nonempty paths; $c_1$ does not
enter their tests. Corollary~\ref{oa:row} therefore gives the
explicit profile \eqref{oa:explicit-p}, a strict spectral gap and
normalized mixing for every probability kernel in the first species.
No pointwise positive lower bound on the gate or kernel is needed.

For $c_1=2,c_2=1,\alpha_0=d=1$ and a smooth bump $g$ supported in
$(1/4,3/4)$, the choice $\omega=3$ meets the test with $W=1$.
The choice $\omega=1$ is precisely the locked example of
Proposition~\ref{prop:smooth-lock} (take $b$ constant there).
Now $P_{21}=0$ at every window length and an arbitrary trace in the
inactive species persists on the gate-free arc. Thus the traveling gate restores communication away from resonance,
while the receiving scan vanishes at resonance.
For constant $b>0$ in the nonresonant example,
Theorem~\ref{nl:th:application} and \eqref{oa:minorant} give
convergence of every nonzero nonnegative logistic solution to
$b\one$; the resonant example fails this conclusion for $(1,0)$ data.

\begin{proposition}[A heterogeneous periodic gate with an unknown profile]
\label{prop:heterogeneous-gate}
Take $L=T=1$, $(c_1,c_2)=(2,1)$, $(d_1,d_2)=(1,0)$,
and any $L^1$ probability kernel in the active component. Put
$g(y)=1-\cos(2\pi y)$ and, for $0<\varepsilon\le10^{-3}$, set
\begin{equation*}
 B_\varepsilon(t,x)=
 \begin{pmatrix}
 -1+\varepsilon\cos(2\pi x)&1\\
 g(x+3t)&-g(x+3t)-\varepsilon\cos(2\pi x)
 \end{pmatrix}.
\end{equation*}
The two-tree and coefficient gap tests hold with $W=1$.
Consequently there is a simple positive Floquet profile and exponential
principal mixing. The row sums depend on $x$, and no positive Floquet
profile is spatially constant.
\end{proposition}
\begin{proof}
With active index $a=1$, the receiving scan for $1\to2$ is
\[
 \int_0^1g(y+3s+2r)\,dr=1,
\]
whereas the sending scan for $2\to1$ is one. Thus $p_*=q_*=1$,
$H_R=N_R=2$, and the diagonal perturbation changes neither tree.
The extreme row sums are $b_+=\varepsilon$ and
$b_-=-\varepsilon$, so $\Delta_b=2\varepsilon$.
Since $0\le g\le2$, the choice $h_+=2+2\varepsilon$ in
Corollary~\ref{oa:row} gives
\[
 \eta_R=e^{-4-4\varepsilon},\qquad
 N_R\Delta_b=4\varepsilon<e^{-4-4\varepsilon}
                 <-\log(1-e^{-4-4\varepsilon}).
\]
For an elementary uniform check, $4+4\varepsilon<5$ and
$e<3$ give $e^{-4-4\varepsilon}>3^{-5}=1/243>1/250
\ge4\varepsilon$. The upper bound $\eta_R\le e^{-4}<1/2$ is immediate.
The gap is therefore verified before any eigenfunction is constructed.
In particular
\[
 r(V_0)\le e^\varepsilon(1-e^{-4-4\varepsilon})^{1/2}
          <e^{-\varepsilon}\le r(M).
\]
If a positive Floquet profile were $p(t,x)=v(t)$, its first component
equation would read
$v_1'=-v_1+v_2+\varepsilon\cos(2\pi x)v_1-\lambda v_1$.
The mild equation makes $v_1$ continuously differentiable in this case.
Its values at $x=0$ and $x=1/2$ contradict $v_1>0$.
Thus every positive Floquet profile must vary in space.
\end{proof}

\section{Quadratic self-regulation and the variance structure}\label{sec:application}

Principal normalization yields the threshold and critical relaxation.
He--Wang \cite[Theorem~3.1]{HW} prove the related $\mathcal R_0$
dichotomy for vector-bias malaria with time-varying delays in an
almost periodic environment: for $\mathcal R_0<1$, the disease-free
almost periodic solution is globally asymptotically stable; for
$\mathcal R_0>1$, a unique positive almost periodic solution is
globally asymptotically stable for initial data with both infection
components strictly positive.

\subsection{Threshold dynamics and the critical branch}
Let $X=C(\mathbb T_L;\mathbb R^m)$ with component maximum norm.
Let
\[
 A(t)u=C\partial_xu+D_d(K-\Id)u+B(t,x)u
\]
have arbitrary real constant velocities $c_i$, rates $d_i\ge0$,
$L^1$ probability kernels, and a continuous time-space periodic Metzler
matrix $B$. Use the mild transport interpretation above. The conclusions
conditional on a positive profile do not require a graph hypothesis.
Write $U(t,s)$ for
the positive linear evolution and $T$ for a time period. Suppose there
is a continuous positive periodic order unit $p(t)$ with exponent
$\lambda$:
\[
 U(t,s)p(s)=e^{\lambda(t-s)}p(t),\qquad
 p_-:=\min_{t,x,i}p_i(t,x)>0,\quad
 p_+:=\max_{t,x,i}p_i(t,x)<\infty .
\]
In the autonomous nonzero-drift theorem one can take $p(t)=\phi$
and any period $T$. In the time-periodic theorem this order unit is
constructed by its exact period-map separation criterion.
No eigenfunction is assumed in applying the extinction conclusion
for a strictly negative actual uniform linear growth bound.

Consider the explicitly cooperative population system
\begin{equation}\label{nl:eq:nonlinear}
 u_t=A(t)u-u\odot u,\qquad
 (u\odot u)_i=u_i^2 .
\end{equation}
The componentwise quadratic term gives a cooperative KPP self-regulation
class. The theorem below establishes its threshold and attraction
properties directly in the mild transport setting.

\begin{theorem}\label{nl:th:application}
Every nonnegative initial function has a unique nonnegative global
mild solution of \eqref{nl:eq:nonlinear}.

If the actual uniform linear growth bound is strictly negative,
every solution tends uniformly to zero, without needing a positive
Floquet function.

If a positive $p$ exists with $\lambda=0$, every nonnegative solution
has $\|u(t)\|=O((1+t)^{-1})$. If $p$ exists with $\lambda<0$, the same
comparison gives exponential decay.

If $p$ exists with $\lambda>0$, there is exactly one strictly positive
continuous $T$-periodic solution $u^*(t)$. Every strictly positive
initial function converges uniformly to it at an exponential rate
depending on its positive floor and cap. If, in addition, the
linear evolution has a finite-time integral minorant
\begin{equation}\label{nl:eq:minorant}
 U(s+H,s)f\ge\eta\one
       \sum_i\int_0^L f_i(x)\,dx\quad(s\in\R,\ f\ge0)
\end{equation}
for some $H,\eta>0$, every nonzero nonnegative initial function
has the same convergence conclusion. No convergence rate uniform
through vanishing initial mass is asserted.

For $\lambda>0$, $u^*$ is also the only entire solution among those
with a finite uniform cap and a joint space-time positive floor.
In the autonomous case this positive periodic solution is independent
of time and is the unique strictly positive steady state.
\end{theorem}

Subsection~\ref{app:nonlinear-proof} proves
Theorem~\ref{nl:th:application} by mild comparison, contraction in an
invariant order interval and minorization for nonzero seeds.

Assume that the positive Floquet profile has exponent zero,
$\ell_t(p(t))=1$, and normalized exponential mixing. In the operator
notation below, $p(t)\ell_s$ means $f\mapsto p(t)\ell_s(f)$:
\begin{equation}\label{nl:eq:critical-mix}
 \|Q(t,s)\|\le C e^{-\kappa(t-s)},\qquad
 Q(t,s)=U(t,s)-p(t)\ell_s,\quad t\ge s.
\end{equation}
The functionals are positive and weak-star continuous and periodic, as
constructed above. Put
\begin{equation*}
 a_*(t)=\ell_t(p(t)\odot p(t)),\qquad
 \bar a_* =\frac1T\int_0^T a_*(t)\,dt.
\end{equation*}
They satisfy $p_-\le a_*(t)\le p_+$, so $\bar a_*>0$.
Rescaling $p$ by a positive constant and $\ell_t$ by its reciprocal
rescales $\bar a_*$ by that constant, leaving $p/\bar a_*$ unchanged.
The coefficient is invariant under simultaneous spatial translation by
the conjugacy already proved.

\begin{lemma}[Projection identities for mild solutions]
\label{nl:lem:mild-projection}
Under \eqref{nl:eq:critical-mix}, suppose that $v$ satisfies
\[
 v(t)=U(t,s)v(s)+\int_s^tU(t,r)F(r)\,dr,
 \qquad F\in C([s,t];X).
\]
Then $a(t)=\ell_t(v(t))$ is continuously differentiable and
\begin{align*}
 a'(t)&=\ell_t(F(t)),\\
 v(t)-a(t)p(t)
 &=Q(t,s)v(s)+\int_s^tQ(t,r)F(r)\,dr.
\end{align*}
No time derivative of the dual functionals, or density for their
representing measures, is needed.
\end{lemma}
\begin{proof}
Positivity and $\ell_r(p(r))=1$ give
$\|\ell_r\|=\ell_r(\one)\le p_-^{-1}$ uniformly.
A fixed bounded functional commutes with a strong Bochner integral.
Apply $\ell_t$ to the mild formula and use
$\ell_tU(t,r)=\ell_r$ at zero exponent to obtain
\[
 a(t)=a(s)+\int_s^t\ell_r(F(r))\,dr.
\]
The integrand is continuous: split its increment into the strong
increment of $F$ and the weak-star increment of $\ell_r$ on the fixed
value of $F$. The uniform functional bound controls the first part.
The fundamental theorem of calculus proves differentiability.
Subtracting the identity multiplied by $p(t)$ from the mild formula
proves the transverse formula.
\end{proof}

\begin{theorem}[Critical amplitude and a small positive growth shift]
\label{nl:th:critical-profile}
Suppose \eqref{nl:eq:critical-mix} holds. Every strictly positive initial
function for \eqref{nl:eq:nonlinear} satisfies
\begin{equation}\label{nl:eq:critical-profile}
 \left\|t u(t)-\frac{p(t)}{\bar a_*}\right\|
                =O\!\left(\frac{\log(2+t)}{1+t}\right).
\end{equation}
The constant may depend on the initial function. If the finite-time
integral minorant \eqref{nl:eq:minorant} also holds, the same assertion
applies to every nonzero nonnegative initial function.
For the equation with $A(t)$ replaced by $A(t)+\varepsilon\Id$,
let $u_\varepsilon^*$ be its positive periodic state. As
$\varepsilon\downarrow0$,
\begin{equation}\label{nl:eq:small-branch}
 \sup_t\left\|u_\varepsilon^*(t)
                  -\frac{\varepsilon}{\bar a_*}p(t)\right\|
                  =O(\varepsilon^2).
\end{equation}
In the autonomous normalization of Theorem~\ref{th:principal},
$p=\phi$ and
$\bar a_* =\int_0^L\sum_i\psi_i\phi_i^2\,dx$.
These conclusions apply to the autonomous mixing class, to the
all-active periodic class with spectral separation, and to the
two-tree class with spectral separation.
\end{theorem}
\begin{proof}
The minorant and the nonlinear lower linear comparison in the proof of
Theorem~\ref{nl:th:application} make the solution a strict order unit at
a finite time. If it is initially strictly positive, start immediately.
The scalar multiples of $p$ with amplitudes solving
$a'=-p_+a^2$ and $a'=-p_-a^2$ are respectively lower and upper
solutions. After a finite shift this gives
\begin{equation}\label{nl:eq:critical-two-sided}
 \frac{c}{1+t}p(t)\le u(t)\le\frac{C}{1+t}p(t).
\end{equation}
Write $a(t)=\ell_t(u(t))$ and $w(t)=u(t)-a(t)p(t)$.
Lemma~\ref{nl:lem:mild-projection}, with $F=-u\odot u$, gives
\[
 a'=-\ell_t(u\odot u),\qquad
 w(t)=Q(t,s)u(s)-\int_s^tQ(t,r)(u(r)\odot u(r))\,dr.
\]
By \eqref{nl:eq:critical-two-sided} and exponential mixing,
$\|w(t)\|=O((1+t)^{-2})$. To check the convolution bound, split
at $r=t/2$: the early integral is exponentially small and the late
one is bounded by a constant times $(1+t)^{-2}$.
Consequently
\[
 a'=-a_*(t)a^2+O((1+t)^{-3}),\qquad
 a(t)\ge c/(1+t).
\]
Divide by $a^2$ and integrate. Periodicity gives
\[
 \frac1{a(t)}=\bar a_* t+O(\log(2+t)),\qquad
 a(t)=\frac1{\bar a_*t}+O\!\left(\frac{\log t}{t^2}\right).
\]
Combining this with the transverse estimate proves
\eqref{nl:eq:critical-profile}.

For the positive shift, the existing existence and uniqueness theorem
applies with principal exponent $\varepsilon$ and the same $p$.
The periodic lower and upper functions
$\varepsilon p/p_+$ and $\varepsilon p/p_-$ trap its unique state;
this also follows by attraction from data in that invariant interval.
Thus $u_\varepsilon^*$ and $a_\varepsilon=\ell_t(u_\varepsilon^*)$
are bounded above and below by positive constants times $\varepsilon$.
For $0<\varepsilon<\kappa/2$, the transverse formula over $n$ periods,
using $U_\varepsilon(t,s)=e^{\varepsilon(t-s)}U(t,s)$, has an initial
term tending to zero as $n\to\infty$. Its integral is bounded by
$C\varepsilon^2\int_0^\infty e^{-(\kappa-\varepsilon)r}dr$.
Hence $\|w_\varepsilon(t)\|=O(\varepsilon^2)$ uniformly in phase.
The amplitude equation is
\[
 a_\varepsilon'=\varepsilon a_\varepsilon
                   -a_*(t)a_\varepsilon^2+O(\varepsilon^3).
\]
It first gives $a_\varepsilon(t)=a_\varepsilon(0)+O(\varepsilon^2)$
over a period. Integrating it over that period and using
$a_\varepsilon(0)\asymp\varepsilon$ yields
$a_\varepsilon(0)=\varepsilon/\bar a_*+O(\varepsilon^2)$.
The transverse estimate proves \eqref{nl:eq:small-branch}.
\end{proof}

\subsection{The cubic amplitude and the logarithmic correction}

Keep the zero-exponent normalization of
Theorem~\ref{nl:th:critical-profile}. Define the transverse response
\begin{equation}\label{nl:eq:response}
 r_*(t)=-\int_0^\infty Q(t,t-s)
              [p(t-s)\odot p(t-s)]\,ds,
 \qquad b_*(t)=-2\ell_t(p(t)\odot r_*(t)).
\end{equation}
The integral converges in $X$. Both functions are continuous and
$T$-periodic, and $\ell_t(r_*(t))=0$. Put
\[
 \mathfrak a=\bar a_*,\qquad
 \mathfrak b=\frac1T\int_0^T b_*(t)\,dt,
\]
and let $J_*$ be the unique mean-zero $T$-periodic primitive of
$a_* -\mathfrak a$. For the nonzero trajectories in
Theorem~\ref{nl:th:critical-profile}, the coefficient identities of
Theorem~\ref{nl:th:normal-form} read
\[
 a'=-a_*a^2+b_*a^3+O(a^4),\qquad
 (a^{-1})'=a_*-b_*a+O(a^2),\quad a=\ell_t(u).
\]
Thus the reciprocal amplitude has logarithmic coefficient
$-\mathfrak b/\mathfrak a$, while the density correction has coefficient
$\mathfrak b p/\mathfrak a^3$, as proved in
Theorem~\ref{nl:th:normal-form}.

\begin{lemma}[The centered response and its characteristic domain]
\label{nl:lem:response-domain}
Let $\mathcal E$ and $\mathscr G$ be the zero-exponent evolution
semigroup and generator, and set $q=p^2-a_*p$. The response
\eqref{nl:eq:response} is continuous, periodic, and satisfies
\[
 \ell_t(r_*(t))=0,\qquad
 r_*\in D(\mathscr G),\qquad \mathscr Gr_*=q.
\]
It is the unique solution of this equation that is centered at every
time phase. Define $\widehat{\mathcal E}(s)f=p^{-1}\mathcal E(s)(pf)$
and let $\widehat{\mathscr G}$ be its generator. Then $h=r_*/p$
belongs to $D(\widehat{\mathscr G})$ and
$\widehat{\mathscr G}h=p-a_*\one$.
\end{lemma}
\begin{proof}
The closed subspace
$Z=\{f\in Y:\ell_t(f(t))=0\text{ for every }t\}$ is invariant
under $\mathcal E$. On $Z$, \eqref{nl:eq:critical-mix} gives
$\|\mathcal E(s)f\|\le Ce^{-\kappa s}\|f\|$.
The function $q$ belongs to $Z$, and dual covariance gives
$\mathcal E(s)q(t)=Q(t,t-s)p(t-s)^2$. Thus
\[
 r_*=-\int_0^\infty\mathcal E(s)q\,ds\quad\hbox{in }Y.
\]
Strong continuity and an exponential tail prove continuity and
centering. The generator difference quotient is
\[
 \frac{\mathcal E(\tau)r_*-r_*}{\tau}
 =\frac1\tau\int_0^\tau\mathcal E(s)q\,ds\longrightarrow q
 \quad\hbox{in }Y.
\]
If two centered domain solutions existed, their difference would be
fixed by $\mathcal E$ and also decay exponentially on $Z$, hence vanish.
Multiplication by the positive continuous order unit $p$ is boundedly
invertible on $Y$. The generator of its conjugated semigroup has domain
$\{f:pf\in D(\mathscr G)\}$, which proves the last assertion without
requiring separate time or space derivatives.
\end{proof}

\begin{lemma}[Cubic dissipation and stationary variance]
\label{nl:lem:dissipation}
Under the zero-exponent mixing normalization above, on $Y=C(\T_T\times\T_L;\R^m)$ use the evolution semigroup
$(\mathcal E(s)f)(t)=U(t,t-s)f(t-s)$ and its generator $\mathscr G$,
as in Proposition~\ref{tp:th:general}. Set
\[
 g_i(t,x)=p_i(t,x)-a_*(t),\quad h=r_*/p,\quad
 \widehat{\mathcal E}(s)f=p^{-1}\mathcal E(s)(pf),\qquad
 \mu(f)=\frac1T\int_0^T\ell_t(p(t)f(t))\,dt.
\]
Products, squares and division by $p$ are componentwise;
$\widehat{\mathscr G}$ denotes the generator of
$\widehat{\mathcal E}$. Then $r_*\in D(\mathscr G)$,
$\mathscr Gr_*=p^2-a_*p$, and
\begin{align}
 \mathfrak b&=-2\mu(h\widehat{\mathscr G}h)
 =\lim_{s\downarrow0}\frac1s\mu\!\left(
 \widehat{\mathcal E}(s)(h^2)
          -(\widehat{\mathcal E}(s)h)^2\right)\ge0,
                       \label{nl:eq:cubic-dissipation}\\
 \mathfrak b&=2\int_0^\infty
                \mu\bigl(g\widehat{\mathcal E}(s)g\bigr)\,ds.
                       \label{nl:eq:cubic-correlation}
\end{align}
The last integral converges absolutely. The positive unital semigroup $\widehat{\mathcal E}$ has a stationary
Markov realization $(\Xi_s)$ with invariant law $\mu$ on
$(\R/T\mathbb Z)\times\T_L\times\{1,\ldots,m\}$. It satisfies
\begin{equation}\label{nl:eq:cubic-variance}
 \mathfrak b=\lim_{S\to\infty}\frac1S
 \operatorname{Var}_{\mu}\!\left(\int_0^Sg(\Xi_s)\,ds\right).
\end{equation}
In the autonomous class of
Theorem~\ref{th:principal}, with $s=0$ and every $d_i>0$,
\begin{equation*}
 \mathfrak b=0\quad\Longleftrightarrow\quad
 \phi_i(x)\equiv k>0\quad\hbox{for every }i,x.
\end{equation*}
This equality characterization is restricted to that autonomous class.
\end{lemma}
\begin{proof}
Lemma~\ref{nl:lem:response-domain} gives the domain statements and
$\widehat{\mathscr G}h=g$. Continuity of $b_*$ follows from the
continuous response, weak-star continuity of $\ell_t$, and its uniform
norm bound in Lemma~\ref{nl:lem:mild-projection}.
Principal conjugation is positive and unital. The functional $\mu$
is a probability, invariant by dual covariance and a periodic change
of variable. Since $\ell_t(r_*)=0$,
$-2\mu(hg)=\mathfrak b$.
Jensen's inequality gives
$(\widehat{\mathcal E}(s)h)^2\le\widehat{\mathcal E}(s)(h^2)$;
invariance and norm differentiation at zero prove
\eqref{nl:eq:cubic-dissipation}. No derivative of the measure-valued
dual or generator-domain assertion for $h^2$ is used here.
Moreover $\ell_t(p(t)g(t))=0$ at every phase, and therefore
\[
 (\widehat{\mathcal E}(s)g)(t)
 =p(t)^{-1}Q(t,t-s)[p(t-s)g(t-s)].
\]
Principal mixing bounds its norm by $Ce^{-\kappa s}$ uniformly in $t$.
The centering is phasewise; a function of time alone need not decay
under this evolution semigroup.
Dividing the stable integral by $p$ gives
$h=-\int_0^\infty\widehat{\mathcal E}(s)g\,ds$.
Boundedness of $g,\mu$ justifies the absolutely convergent pairing
and proves \eqref{nl:eq:cubic-correlation}; no pointwise sign of its
integrand is needed.

$\widehat{\mathcal E}$ is a conservative Feller
semigroup on the displayed compact metric space: it is strongly
continuous on its continuous functions, positive and unital.
The Feller realization theorem
\cite[Chapter~4, Theorem~2.7]{EK} supplies a process with right-continuous
paths and left limits. Starting it with the invariant probability $\mu$
makes it stationary. Since $\mu(g)=0$, boundedness of $g$ justifies
Fubini on each finite square and the Markov property gives
\[
 \frac1S\operatorname{Var}_{\mu}\!\left(\int_0^Sg(\Xi_s)\,ds\right)
 =2\int_0^S(1-s/S)\mu\bigl(g\widehat{\mathcal E}(s)g\bigr)\,ds.
\]
The absolute integrability already proved gives
\eqref{nl:eq:cubic-variance} by dominated convergence.
This is the classical integrated-covariance identity, here identifying
the nonlinear correction coefficient; it requires no reversibility or
mixing of functions depending only on the periodic phase.

For the autonomous assertion write $p=\phi$,
$\pi(f)=\sum_i\int_0^L\psi_i\phi_i f_i$ and
$\widehat Af=\phi^{-1}A(\phi f)$.
The response is time independent, and the same norm difference quotient
for the autonomous semigroup puts $r_*$ in $D(A)=C^1$.
Thus $h\in C^1$, $h^2\in D(\widehat A)$ and invariance gives
$\mathfrak b=\pi(\widehat A(h^2)-2h\widehat Ah)$.
Expanding the transport derivative by its usual product rule cancels
that part, leaving exactly
\begin{equation*}
 \begin{split}
 \mathfrak b={}&\sum_i d_i\int_0^L\!\int_\R
 \psi_i(x)\phi_i(x-z)J_i(z)
       [h_i(x-z)-h_i(x)]^2\,dz\,dx\\
 &+\sum_{i\ne j}\int_0^L\psi_i(x)B_{ij}(x)\phi_j(x)
                         [h_j(x)-h_i(x)]^2\,dx .
 \end{split}
\end{equation*}
All integrals are finite and nonnegative. If their sum is zero,
strict positivity of $\psi,\phi$ and $d_i$ implies
$K_i h_i=h_i$ for every $x$: the squared jump difference has zero
integral at every $x$ by continuity of its convolution expression.
Iterate this identity to a finite convolution power with a positive
circle minorant, furnished by Lemma~\ref{le:convolution}.
Its oscillation contraction forces $h_i$ to be constant in space.
For every edge with $\int B_{ij}>0$, the second line then forces
the two constants to agree. Union-graph connectivity makes them all
equal. Since $\pi(h)=\ell(r_*)=0$ and $\pi(\one)=1$, they are zero.
The response equation $\widehat Ah=\phi-a_*\one$ therefore gives
the common constant profile. Conversely such a profile has
$p^2=a_*p$, hence $r_*=0$ and $\mathfrak b=0$.
\end{proof}

\begin{theorem}[Second-order critical and positive-state expansions]
\label{nl:th:normal-form}
Under the hypotheses and initial-data alternatives of
Theorem~\ref{nl:th:critical-profile}, set $a(t)=\ell_t(u(t))$.
As $t\to\infty$,
\begin{align}
 u(t)-a(t)p(t)&=a(t)^2r_*(t)+O(a(t)^3),\label{nl:eq:transverse-second}\\
 a'&=-a_*(t)a^2+b_*(t)a^3+O(a^4),\label{nl:eq:cubic}\\
 \frac1{a(t)}&=\mathfrak a t-
       \frac{\mathfrak b}{\mathfrak a}\log t+J_*(t)+C(u_0)
                  +O(\log t/t).\label{nl:eq:inverse-second}
\end{align}
In particular the full profile expansion is
\begin{equation}\label{nl:eq:profile-second}
 u(t)=\frac{p(t)}{\mathfrak a t}
 +\frac1{t^2}\left\{
 \frac{\mathfrak b}{\mathfrak a^3}(\log t)p(t)
 +\frac{r_*(t)-[J_*(t)+C(u_0)]p(t)}{\mathfrak a^2}
 \right\}+O((\log t)^2/t^3).
\end{equation}
For the positive periodic state of $A(t)+\varepsilon\Id$,
uniformly in $t$,
\begin{equation}\label{nl:eq:branch-second}
 u_\varepsilon^*(t)=\frac{\varepsilon}{\mathfrak a}p(t)
 +\varepsilon^2\left\{
 \frac{\mathfrak b}{\mathfrak a^3}p(t)
 +\frac{r_*(t)-J_*(t)p(t)}{\mathfrak a^2}
 \right\}+O(\varepsilon^3).
\end{equation}
In particular its period-averaged principal amplitude satisfies
\begin{equation}\label{nl:eq:mean-branch}
 \frac1T\int_0^T\ell_t(u_\varepsilon^*(t))dt
 =\frac{\varepsilon}{\mathfrak a}
  +\frac{\mathfrak b}{\mathfrak a^3}\varepsilon^2
  +O(\varepsilon^3).
\end{equation}
The averaged cubic coefficient satisfies $\mathfrak b\ge0$.
All errors are uniform in the periodic phase. Critical-trajectory
constants may depend on its initial function.
\end{theorem}
\begin{proof}
The previous theorem gives $a\asymp(1+t)^{-1}$,
$w=u-ap=O((1+t)^{-2})$, and $|a'|\le C(1+t)^{-2}$.
Consequently $u^2=a^2p^2+O((1+t)^{-3})$, with componentwise
products understood. Compare the exact transverse Duhamel formula
with $a(t)^2r_*(t)$. In its recent part $t/2\le s\le t$,
\[
 |a(s)^2-a(t)^2|\le C(t-s)(1+t)^{-3}.
\]
Integrating this bound against $Ce^{-\kappa(t-s)}$, as well as
the cubic nonlinear error, gives $O((1+t)^{-3})$. The part
$s\le t/2$, the initial term and the omitted tail of $r_*$
are exponentially small. This proves
\eqref{nl:eq:transverse-second}, since $a\ge c/(1+t)$.
The exact mild amplitude identity remains
$a'=-\ell_t(u^2)$. Substitute the transverse expansion:
the cross term is $-2a^3\ell_t(pr_*)$, while the other terms are
$O(a^4)$. This proves \eqref{nl:eq:cubic} without differentiating
the dual functionals or $r_*$.

Dividing by $a^2$ and using the first-order expansion for $a$ gives
\[
 (1/a)'=a_*(t)-\frac{b_*(t)}{\mathfrak a t}
                          +O(\log t/t^2).
\]
A zero-mean continuous periodic function has a bounded primitive.
Integration by parts therefore shows that
$\int^t[b_*(s)-\mathfrak b]s^{-1}ds$ has a finite limit with
$O(t^{-1})$ tail. The remaining error has an integrable tail of
order $O(\log t/t)$. Integrating the displayed identity proves
\eqref{nl:eq:inverse-second} with one finite constant $C(u_0)$.
Expanding its reciprocal, whose correction is $O(\log t)$,
and substituting into \eqref{nl:eq:transverse-second} proves
\eqref{nl:eq:profile-second}.

For the periodic positive state write $a_\varepsilon=\ell_t(u_\varepsilon^*)$
and $w_\varepsilon=u_\varepsilon^*-a_\varepsilon p$.
The preceding theorem gives $a_\varepsilon\asymp\varepsilon$,
$w_\varepsilon=O(\varepsilon^2)$ and
$|a_\varepsilon'|\le C\varepsilon^2$ at every phase.
The transverse formula over infinitely many past periods is
\[
 w_\varepsilon(t)=-\int_0^\infty e^{\varepsilon s}
                  Q(t,t-s)[u_\varepsilon^*(t-s)^2]ds.
\]
For $0<\varepsilon<\kappa/2$, replacing the squared amplitude at
$t-s$ by that at $t$ costs at most $C\varepsilon^3s$;
replacing $e^{\varepsilon s}$ by one costs at most
$\varepsilon s e^{\varepsilon s}$ times $C\varepsilon^2$.
The cross terms are $O(\varepsilon^3)$. Exponential integration
therefore gives, uniformly in phase,
\[
 w_\varepsilon=a_\varepsilon^2r_*+O(\varepsilon^3),\qquad
 a_\varepsilon'=\varepsilon a_\varepsilon
      -a_*a_\varepsilon^2+b_*a_\varepsilon^3+O(\varepsilon^4).
\]
Using $a_\varepsilon=\varepsilon/\mathfrak a+O(\varepsilon^2)$
in the last identity first gives
$a_\varepsilon'=-\varepsilon^2J_*'/\mathfrak a^2+O(\varepsilon^3)$.
For this paragraph write $\overline f=T^{-1}\int_0^T f(t)dt$.
The first-order bound makes
$c_\varepsilon=\varepsilon^{-2}(\overline{a_\varepsilon}
 -\varepsilon/\mathfrak a)$ bounded. Integrating the last derivative
estimate over a period gives
\[
 a_\varepsilon(t)=\varepsilon/\mathfrak a
 +\varepsilon^2[c_\varepsilon-J_*(t)/\mathfrak a^2]
 +O(\varepsilon^3).
\]
Put $k_\varepsilon=c_\varepsilon-J_*/\mathfrak a^2$ and average
the amplitude equation. The order-$\varepsilon^2$ mean vanishes, and
\[
 \begin{split}
 0&=\varepsilon^3\left\{
 \overline{k_\varepsilon}
 -\frac{2}{\mathfrak a}\overline{a_*k_\varepsilon}
 +\frac{\mathfrak b}{\mathfrak a^3}\right\}+O(\varepsilon^4)\\
 &=\varepsilon^3\left(-c_\varepsilon
                 +\frac{\mathfrak b}{\mathfrak a^3}\right)
                         +O(\varepsilon^4).
 \end{split}
\]
Here $\overline{J_*}=0$ and
$\overline{a_*J_*}=\overline{(\mathfrak a+J_*')J_*}=0$.
Thus $c_\varepsilon=\mathfrak b/\mathfrak a^3+O(\varepsilon)$.
Combining the amplitude and transverse estimates proves
\eqref{nl:eq:branch-second}. Pairing with $\ell_t$, using
$\ell_t(r_*)=0$, $\ell_t(p)=1$ and $\overline{J_*}=0$, and
averaging the uniform remainder gives \eqref{nl:eq:mean-branch}.
\end{proof}

The transverse $r_*$ and phase $J_*$ contribute to the full vector
term; $C(u_0)$ records its initial state. If $\mathfrak b=0$, the
displayed logarithmic term vanishes. Lemma~\ref{nl:lem:dissipation}
classifies that equality in the autonomous all-active class; no next
nonzero order is identified in the general periodic or partial class.

\begin{proposition}
\label{nl:prop:log-sharp}
On the unit circle take a smooth strictly positive probability kernel
whose periodization is one, let $P$ be averaging, and put
$p(x)=1+\eta\cos(2\pi x)$, $0<|\eta|<1$. The scalar operator
\[
 A=P-p^{-1}\Id=(P-\Id)+(1-p^{-1})\Id
\]
has exponent zero, right function $p$, dual density
$p/\int p^2$ and exponential principal mixing. Its cubic coefficient is
\[
 \mathfrak a=\frac{\int p^3}{\int p^2},\qquad
 \mathfrak b=\frac{2}{\int p^2}\int p^3(p-\mathfrak a)^2>0.
\]
Consequently every nonzero nonnegative solution of $u_t=Au-u^2$
satisfies
\[
 \frac{t}{\log t}\left[tu(t)-\frac p{\mathfrak a}\right]
       \longrightarrow \frac{\mathfrak b}{\mathfrak a^3}p
 \quad\hbox{in }X.
\]
For $\eta=1/2$ the exact values are
$\mathfrak a=11/9$ and $\mathfrak b=893/5832$.
\end{proposition}
\begin{proof}
Convolve a Gaussian with the uniform density on $[-1/2,1/2]$ to
realize $P$. The identity $Ap=0$, symmetry under integration and
$\int p=1$ give the claimed principal data: directly
$\int pAf=\int (Ap)f=0$. After conjugating by $p$, the generator
is
\[
 Lf(x)=\frac1{p(x)}\left(\int_0^1p(y)f(y)dy-f(x)\right).
\]
It is Markov. With $q_-=1/\max p$, $q_+=1/\min p$, its one-jump
positive term gives, for every $t>0$ and $f\ge0$,
$e^{tL}f\ge t q_- e^{-tq_+}\int p f$.
This follows by retaining a single $qP_p$ between two diagonal
retentions $e^{-q s}$ and integrating its jump time from zero to $t$.
Markov oscillation contraction proves exponential mixing; undoing
the conjugacy also supplies a full-support minorant for $e^{tA}$.
The scalar, all-active periodic class allows its zero velocity.

Write $I_j=\int_0^1p^j$. The unique transverse response solves
$Ar_*=p^2-\mathfrak a p$ and $\int pr_*=0$; explicitly
\[
 r_*=mp-p^3+\mathfrak a p^2,
 \qquad m=(I_4-\mathfrak a I_3)/I_2.
\]
Integration verifies $\int r_*=m$, so this is the exact equation,
not just its projection. Hence
$-2\int p^2r_*/I_2
=2(I_5-2\mathfrak a I_4+\mathfrak a^2I_3)/I_2$,
which is the asserted positive integral. Positivity is strict because
$p$ is positive and nonconstant. For $\eta=1/2$,
$(I_2,I_3,I_4,I_5)=(9/8,11/8,227/128,303/128)$,
giving the two stated rational values. The profile conclusion now
follows from \eqref{nl:eq:profile-second}; all bounded periodic terms
vanish after its indicated rescaling.
\end{proof}

\paragraph{A constant two-species check.}
On the unit circle, take any nonzero component velocities, positive
jump rates and probability kernels, and
$B=\left(\begin{smallmatrix}-2&1\\2&-1\end{smallmatrix}\right)$.
Then $p=(1,2)^T$ and
$\ell(f)=\tfrac13\int_0^1(f_1+f_2)$ give the zero-growth
mixing normalization. Directly
\[
 \mathfrak a=5/3,\quad
 p^2-\mathfrak a p=(-2/3,2/3)^T,\quad
 r_*=(2/9,-2/9)^T,\quad \mathfrak b=4/27.
\]
Indeed $Br_*=p^2-\mathfrak a p$ and $\ell(r_*)=0$;
spatially constant functions have no transport or dispersal terms.
The unique small positive state therefore has expansion
\[
 u_\varepsilon^*
 =\varepsilon\binom{3/5}{6/5}
  +\varepsilon^2\binom{14/125}{-2/125}+O(\varepsilon^3).
\]
Thus species variation alone makes the logarithmic coefficient positive,
whereas the full vector's second correction has opposite signs.

\begin{remark}\label{rem:nonlinear-threshold-scope}
The principal growth rate measures invasion at low density, while the
profile records spatial and species composition. At critical growth,
mixing selects the $t^{-1}$ composition and its logarithmic relaxation
correction. Just above criticality, the small positive branch quantifies
the abundance sustained by quadratic regulation. Its averaged principal
amplitude has a nonnegative second correction, although individual
component corrections can have different signs. These interpretations
require the proved principal-profile and mixing hypotheses: a zero
actual growth bound alone does not yield the critical decay comparison,
and a positive bound alone does not ensure a positive periodic state.
The model and conclusions concern a periodic habitat.
\end{remark}

\section{Examples and scope of the conclusions}\label{sec:examples-scope}

The examples below probe what happens when the velocity, communication,
or spectral-separation conditions are weakened. They also show where the
sufficient criteria leave room for other communication mechanisms.

\subsection{Examples and the loss of a transport direction}\label{sec:examples}
Propositions~\ref{prop:smooth-lock}--\ref{prop:dual-lock} distinguish
the two scan failures and quantify the delay near resonance.
Proposition~\ref{prop:drift-crossing} records a drift-induced threshold
crossing; Appendices~\ref{app:stationary}--\ref{sec:floquet} give the
stationary and resonant spectral obstructions.

\paragraph{Disjoint gates and kernels with different tail behavior.}
Take $m=2$, $L=1$, $c_1=1$, $c_2=-2$, $d_1=1$ and any $d_2\ge0$.
Let $B_{12},B_{21}$ be nonzero smooth nonnegative bumps supported,
respectively, in $(0,1/4)$ and $(1/2,3/4)$, and take arbitrary real
continuous periodic diagonal entries. Every pointwise matrix is reducible,
but the union graph is strongly connected. For example, one may use
\[
 J_1(h)=\frac{p-1}{2}(1+|h|)^{-p},\quad 1<p\le2,
 \qquad J_2(h)=e^{-h}\one_{\{h\ge0\}}.
\]
The first kernel has infinite first absolute moment and the second is
asymmetric. All conclusions of Theorems~\ref{th:principal}--\ref{th:perturb}
apply. In particular the active first species suffices even when $d_2=0$.

When all $d_i=0$, mixing is not automatic in the class of
Theorem~\ref{th:principal}. For instance $m=1$, $B=0$ and $c_1\ne0$
give pure periodic translation, which does not converge to its mean for
a nonconstant initial function. This does not exclude mixing by switching
in some zero-jump systems with different velocities.

\subsubsection{Smooth periodic gates and an inactive component}
Frozen autonomous mixing does not imply periodic mixing, even with
smooth coefficients and a positive Floquet function.

\begin{proposition}\label{prop:smooth-lock}
On $\T=\R/\mathbb Z$, choose a nonzero $g\in C^\infty(\T)$ with
$g\ge0$ and support in $(1/4,3/4)$ modulo one. Put
$E=(-1/8,1/8)$ modulo one. Let $K$ be convolution by the strictly
positive smooth probability density $\pi^{-1/2}e^{-h^2}$, and fix $b\ge0$.
Consider
\begin{equation}\label{eq:smooth-lock}
 \begin{split}
 (u_1)_t&=2(u_1)_x+(K-\Id)u_1+(b-1)u_1+u_2,\\
 (u_2)_t&=(u_2)_x+g(x+t)u_1+(b-g(x+t))u_2.
 \end{split}
\end{equation}
For each frozen time, its autonomous generator has velocities $2,1$,
one positive jump rate, and a strongly connected spatial union graph.
It satisfies Theorems~\ref{th:principal} and \ref{th:mix}, with
principal value $b$ and principal function $\one$.

The true one-period maps $M=U(1,0)$ and $V_0$ with the jump term deleted
instead satisfy
\begin{equation}\label{eq:smooth-radii}
 r(M)=r(V_0)=e^b,
 \qquad U(t,s)\one=e^{b(t-s)}\one.
\end{equation}
The normalized powers $e^{-bn}M^n$ cannot converge strongly to an
operator of finite rank. There is no positive fixed-time integral
minorant of the form \eqref{nl:eq:minorant}.
Thus the all-positive jump-rate assumption cannot be dropped from
Theorem~\ref{tp:th:main} merely by retaining one active species.

For $b>0$, the same example with the losses $-u_i^2$ added has the
positive periodic solution $b\one$. A nonzero nonnegative initial
function need not converge to it, although every strictly positive
initial function does.
\end{proposition}
\begin{proof}
All coefficient and kernel assumptions are literal: assign the displayed
Gaussian to both kernels, take $c_1=2,c_2=1$, $d_1=1,d_2=0$, and
\[
 B(t,x)=b\Id+
 \begin{pmatrix}-1&1\\g(x+t)&-g(x+t)\end{pmatrix}.
\]
At every frozen time the two off-diagonal integrals are positive.
The frozen domain is exactly $C^1(\T;\R^2)$, and
$A(t)\one=b\one$ identifies the frozen principal value and function.
The true evolution is the unique positive mild transport evolution on
$C(\T;\R^2)$. Its row sums are $b$, so uniqueness and comparison give
the identity on constants in \eqref{eq:smooth-radii} and
$\norm{U(t,s)}=e^{b(t-s)}$. In particular $p(t)\equiv\one$ is a
strictly positive periodic Floquet function with exponent $b$.

To find the obstruction, follow only the second component's
characteristic. For $f=u(0)$ and fixed $y\in\T$, its mild equation
makes $q_y(t)=u_2(t,y-t)$ continuously differentiable and gives
\[
 q_y'(t)=g(y)u_1(t,y-t)+(b-g(y))q_y(t),\qquad q_y(0)=f_2(y).
\]
This calculation is valid for continuous initial data; no spatial
derivative of $f$ is used. Since $g=0$ on $E$,
\begin{equation}\label{eq:locked-trace}
 u_2(t,y-t)=e^{bt}f_2(y)\qquad(y\in E).
\end{equation}
At integer times this says that restriction of the second component to
any fixed closed subarc $E_0\subset E$ with nonempty interior is unchanged
by $e^{-bn}M^n$. Write that restriction operator as $R_E$; its range is
infinite-dimensional. If the powers converged strongly to a finite-rank
operator $F$, then $R_EF=R_E$, an impossibility. This proves the assertion about
normalized powers, including failure of a rank-one limit.

The no-jump evolution has maximal row sum $b$, so
$\norm{V(t,s)}\le e^{b(t-s)}$. Its second component obeys the same
identity \eqref{eq:locked-trace}. With initial vector $\one$, its norm
at every integer time is at least $e^{bn}$. Periodicity then gives
$r(V_0)=e^b$. The constant eigenvector and the full norm bound give
$r(M)=e^b$, proving \eqref{eq:smooth-radii}.

For any proposed time $H>0$, start with $f=(1,0)$ at time zero and
take $x=y-H$, $y\in E$. Formula~\eqref{eq:locked-trace} gives
$u_2(H,x)=0$, although the total initial integral equals one.
A positive integral minorant is therefore impossible even at this
single initial phase.

Finally add the quadratic losses. On the same characteristics in $E$,
the second component solves $q_y'=bq_y-q_y^2$. Starting with
$f=(1,0)$ makes it zero there for every future time. But $b\one$ is
a positive solution when $b>0$, and the distance from this trajectory
to $b\one$ is at least $b$ at every time. Strictly positive data do
converge to $b\one$ by the positive-profile part of
Theorem~\ref{nl:th:application}; the minorant hypothesis for extending
that conclusion to arbitrary nonzero seeds is absent here.
\end{proof}

The obstruction comes from a time-space gate fixed along one
characteristic, rather than from coefficient roughness. It also explains
why an unqualified substitution of a frozen autonomous mixing theorem
into a periodic evolution is invalid. Both the nonnegative-kernel
definition and the generator domains remain unchanged in this example.

\subsubsection{Delay and characteristic scan obstructions}
\begin{proposition}
\label{prop:near-resonance}
In \eqref{eq:smooth-lock}, replace the second velocity by $1+\delta$,
where $0<|\delta|<1/4$, and keep all coefficients unchanged.
For each such $\delta$ the normalized evolution
$\widehat U_\delta(t,s)=e^{-b(t-s)}U_\delta(t,s)$ has exponential
rank-one convergence. Nevertheless, if
\[
 D_\delta(t)=\sup_{0\le f\le\one}
       \operatorname{osc}(\widehat U_\delta(t,0)f),\qquad
 T_{\mathrm{mix}}(\delta)=\inf\{t:D_\delta(t)\le1/2\},
\]
then
\[
 D_\delta(t)=1\quad(0\le t\le1/(16|\delta|)),\qquad
 T_{\mathrm{mix}}(\delta)\ge1/(16|\delta|).
\]
Here oscillation is over all components and spatial points.
There is also a computable upper bound, with no claim of matching order.
Let $I_g=\int_0^1g$, $W=\lceil1/|\delta|\rceil$,
$h_0=\max(2,\|g\|_\infty)$, $p_0=\min(1,I_g/|\delta|)$,
$\eta_G=\pi^{-1/2}e^{-1/4}$ and
\[
 \epsilon_\delta=\min\{1/2,2\eta_Gp_0W e^{-3h_0W}\}.
\]
Then
\[
 D_\delta(t)\le(1-\epsilon_\delta)^{\lfloor t/(3W)\rfloor},\qquad
 T_{\rm mix}(\delta)\le3W
 \left\lceil\frac{\log2}{-\log(1-\epsilon_\delta)}\right\rceil.
\]
These estimates give an explicit positive exponential rate for each
nonzero $\delta$ and a coarse upper bound $C_g e^{C_g/|\delta|}$
near zero. They do not establish $T_{\rm mix}=O(1/|\delta|)$ or
identify the optimal spectral exponent and projection prefactor.

\end{proposition}
\begin{proof}
The edge $1\to2$ is scanned along receiving velocity $1+\delta$.
Choose the integer-period window $W=\lceil1/|\delta|\rceil$.
Its relative gate coordinate traverses at least one full circle, so
the receiving integral is at least
$|\delta|^{-1}\int_0^1g>0$, uniformly in the starting phase.
The edge $2\to1$ has coefficient one, so its sending scan equals
$W$. These are the two trees with active species one. The row sums
remain $b$. The common-row conclusion of Corollary~\ref{oa:row}
therefore gives spectral separation and normalized mixing.

For a second-component characteristic $x(t)=y-(1+\delta)t$,
the gate coordinate is $y-\delta t$. Choose $y_\pm=\pm1/32$.
For $t\le1/(16|\delta|)$ both coordinates stay in
$(-1/8,1/8)$, where $g=0$. The normalized second-component traces
are consequently their initial values. Choose a continuous
$0\le f\le\one$ with $f_2(y_-)=0$ and $f_2(y_+)=1$.
At each such time the output oscillation is at least one. Since
$\widehat U_\delta$ is positive and preserves $\one$, it cannot
exceed one. This proves the identities and the lower bound.
For the upper estimate apply Lemma~\ref{oa:minorant-return} directly
to the normalized Markov evolution. Its diagonal retention cost can
be chosen as $h_0$. The sending path has coefficient $W$, so its
minimum product including the empty active path is one. The receiving
minimum product is at least $p_0$. On the unit circle the periodized
Gaussian is at least $\eta_G$, by retaining the integer translate
whose argument has modulus at most $1/2$. Thus one active jump
already gives the convolution minorant, and \eqref{oa:constants}
at $m=2,n=1$ gives $\zeta\ge\eta_Gp_0W e^{-3h_0W}$.
Use the probability measure $\nu(f)=\sum_i\int f_i/2$.
Subtracting $\epsilon_\delta\nu(f)\one$ from the Markov map
leaves a positive operator preserving $(1-\epsilon_\delta)\one$.
Its oscillation contraction every $3W$ units is therefore
$1-\epsilon_\delta$, uniformly in phase. Intermediate Markov
evolution cannot increase oscillation. This proves both displayed
upper bounds. Since $p_0=1$ for small $|\delta|$ and
$W\le1/|\delta|+1$, $-\log(1-\epsilon)\ge\epsilon$ gives the
coarse exponential-in-reciprocal-speed estimate.

For a concrete smooth-gate bound, choose $0\le g\le1$
smooth, supported in $(1/4,3/4)$ and equal to one on $[3/8,5/8]$.
At $\delta=1/8$ one has $I_g\ge1/4$, $W=8$, $p_0=1$ and
$h_0=2$. Thus $\epsilon_\delta=16\pi^{-1/2}e^{-1/4-48}$.
This lower bound uses a smooth gate satisfying the coefficient
assumptions.

\end{proof}

The lower delay is the travel time of a locked characteristic; the
upper bound retains one Dyson word on a window of length $W$.

Failure of a tree test alone does not imply failure of mixing.
At $\delta=0$, add any positive jump rate to species two, with its
assigned Gaussian kernel. The row sums are still $b$, so
Corollary~\ref{tp:cor:row} gives mixing. Yet the same receiving
scan for the edge $1\to2$ is zero on the locked arc for every window.
Thus failure of this chosen two-tree certificate is not a converse
nonmixing criterion in the class that requires at least one active
species and allows others to jump. Necessity in the exactly-one-active
class is not inferred from this example.

For the critical choice $b=0$, the same system with quadratic losses
has $t u(t)\to\one$ from every nonzero nonnegative datum when
$\delta\ne0$, by Theorem~\ref{nl:th:critical-profile}: here $p=\one$
and $\bar a_*=1$. At resonance, a datum with second component zero
retains that zero along the locked characteristics even after adding
the quadratic loss. It cannot have this critical profile. The distinction
concerns the actual evolution and the nonzero initial-data class.

\begin{proposition}
\label{prop:dual-lock}
Keep $g,E,K,b$ from Proposition~\ref{prop:smooth-lock} and consider
\begin{equation*}
 \begin{split}
 (u_1)_t&=-2(u_1)_x+(K-\Id)u_1+(b-1)u_1+g(x-t)u_2,\\
 (u_2)_t&=-(u_2)_x+u_1+(b-g(x-t))u_2.
 \end{split}
\end{equation*}
This smooth periodic system has exactly one active species and a
strongly connected spatial union graph at every frozen time.
Its receiving tree $1\to2$ satisfies the test in \eqref{oa:P},
but its sending tree $2\to1$ fails \eqref{oa:Q} for every window.
There is no fixed-time full-support integral minorant. If $M$ is its
one-period map, $e^{-bn}M^n$ cannot converge strongly to any finite-rank
operator.
\end{proposition}
\begin{proof}
The velocities are $-2,-1$, rates $1,0$, and both assigned kernels are
the same smooth Gaussian. The off-diagonal coefficients $g(x-t)$ and
one give the frozen strongly connected graph. The full and no-jump frozen
domains are $C^1(\T;\R^2)$; the mild evolution acts on continuous data.
For every window $W$, the receiving scan of $B_{21}=1$ is $W$.
The sending scan of $B_{12}$ along species two is
$Wg(y-t)$, hence vanishes when $y-t\in E$.

For any nonzero nonnegative smooth $f$ supported in $E$, the exact
solution from $(0,f)$ is
\[
 u_1(t,x)=0,\qquad u_2(t,x)=e^{bt}f(x-t),
\]
as $g(x-t)f(x-t)=0$. Its positive initial integral and identically
zero first component exclude every full-support minorant. At integer
times $e^{-bn}M^n$ fixes all these $(0,f)$. Their span is
infinite-dimensional, whereas any finite-rank strong limit would act
as its identity. This proves failure of that limit.
\end{proof}

Together, Propositions~\ref{prop:smooth-lock} and \ref{prop:dual-lock}
give smooth one-active examples with respectively a sending path alone
and a receiving path alone. Proposition~\ref{prop:near-resonance}
verifies both paths and obtains mixing away from resonance. These
examples exhibit the distinct obstructions addressed by the two scans.

\subsubsection{A heterogeneous example whose principal value changes with drift}

\begin{proposition}
\label{prop:drift-crossing}
On the circle of length $2\pi$, let $P$ denote spatial averaging and set
\[
 A_{c,\eta}=\operatorname{diag}(c,-1)\partial_x+(P-\Id)\Id
 +b\Id+\begin{pmatrix}-1+\eta\cos x&1\\1&-1\end{pmatrix},
 \qquad 1\le c\le2.
\]
Each $P$ is realized by the smooth positive probability density obtained
by convolving the Gaussian with the uniform density on $[-\pi,\pi]$.
For small $\eta$, its principal value satisfies, uniformly on this
velocity interval,
\begin{equation*}
 s(c,\eta)=b+\frac{2\eta^2}{5c^2-2c+13}+O(\eta^4),\qquad
 \partial_cs=-\frac{2(10c-2)\eta^2}{(5c^2-2c+13)^2}+O(\eta^4).
\end{equation*}
In particular, for $b=-3\eta^2/32$ and all sufficiently small
nonzero $\eta$, $s(1,\eta)>0>s(2,\eta)$.
The crossing is unique and its velocity is
$c_{\mathrm{crit}}(\eta)=(3+8\sqrt6)/15+O(\eta^2)$.
For the critically shifted family $A_{c,\eta}-s(c,\eta)\Id$,
normalize $(4\pi)^{-1}\sum_{i=1}^2\int_0^{2\pi}p_i(x)\,dx=1$
and its dual by $\ell(p)=1$. The quadratic self-regulation coefficient satisfies
\begin{equation}\label{eq:drift-cubic}
 \mathfrak b(c,\eta)=
 \frac{c^2+2c+29}{(5c^2-2c+13)^2}\eta^2+O(\eta^4)>0
\end{equation}
for sufficiently small nonzero $\eta$, uniformly on $[1,2]$.
The two leading coefficients at $c=1,2$ are $1/8$ and $37/841$.

\end{proposition}
For the quadratic regulation of Section~\ref{sec:application}, the sign
change is a persistence--extinction transition. After centering at critical
growth, \eqref{eq:drift-cubic} also distinguishes the logarithmic
relaxation corrections. The calculation is in Appendix~\ref{app:drift-calculation}.

\subsection{Scope of the conclusions}

The examples identify the roles of nonzero autonomous velocities,
periodic communication, and strict spectral separation. The two-tree
and row tests are sufficient criteria; they do not classify every
partially dispersive system that mixes. The cubic equality case is
classified only for autonomous all-active systems. All mixing and
self-regulation conclusions concern a compact periodic habitat;
whole-line spectral gaps and propagation speeds require separate
analysis.

\appendix
\section{Technical obstructions and periodic realization}\label{app:technical}

\subsection{A stationary component and concentrated feedback}\label{app:stationary}

\begin{theorem}\label{th:zero}
On $\T=\R/\mathbb Z$ write $Pf=\overline f=\int_\T f(x)\,dx$, and set
\begin{equation*}
 \begin{gathered}
 r(x)=\operatorname{dist}(x,\mathbb Z),\quad
 a(x)=16\sqrt{2r(x)},\\
 \mathcal A\binom uv
  =\binom{-au+\overline u+v}{v'+\overline v-2v+u},
 \quad D(\mathcal A)=C(\T)\times C^1(\T).
 \end{gathered}
\end{equation*}
This is a cooperative dispersal system with $c_1=0$, $c_2=1$,
$d_1=d_2=1$ and both off-diagonal coefficients identically one.
Both averaging operators are realized by a common smooth, even, strictly
positive real-line probability density with every absolute moment finite.
The operator generates a positive $C_0$ semigroup and satisfies
\begin{equation}\label{eq:zero-bound}
 s(\mathcal A)=\omega_0(\mathcal A)=0.
\end{equation}
It has no nonzero nonnegative eigenvector in its domain for any eigenvalue.
Its resolvent is not compact. The coefficient $a$ is continuous and
non-Lipschitz; no conclusion for smoother coefficient classes is asserted.
\end{theorem}

\begin{proof}
\emph{Kernel and realization.}
Let $G(h)=\pi^{-1/2}e^{-h^2}$ and
\[
 J(h)=\int_{-1/2}^{1/2}G(h-y)\,dy.
\]
Convolution with a bounded uniform displacement preserves Gaussian moment
bounds. Thus $J$ is smooth, even, strictly positive, normalized and has
every absolute moment finite. Its periodization equals one: Tonelli and
the tiling by the intervals $[h+n-1/2,h+n+1/2]$ give
$\sum_nJ(h+n)=1$. Its convolution on continuous period-one functions is
therefore $P$. In the form \eqref{eq:A}, the reaction matrix is
\[
 B(x)=\begin{pmatrix}1-a(x)&1\\1&-1\end{pmatrix}.
\]
The diagonal operator $\diag(-a,\partial_x-2)$ generates
$(e^{-ta(x)}u(x),e^{-2t}v(x+t))$ on exactly $C\times C^1$.
The remaining map $(u,v)\mapsto(Pu+v,u+Pv)$ is bounded and positive.
Its Dyson series gives the stated positive semigroup on that unchanged
mixed domain.

\emph{A transport inverse.}
For $\alpha>0$ put
\begin{equation*}
 R_\alpha f(x)=\int_0^\infty e^{-\alpha t}f(x+t)\,dt
     =\frac1{1-e^{-\alpha}}\int_0^1e^{-\alpha t}f(x+t)\,dt.
\end{equation*}
For continuous periodic $f$, a change of variables shows
$(R_\alpha f)'=\alpha R_\alpha f-f$, so
$(\alpha-\partial_x)R_\alpha f=f$ and
\[
 \norm{R_\alpha f}_\infty\le
                 (1-e^{-\alpha})^{-1}\norm f_{L^1(\T)}.
\]
For $A_{22}=\partial_x+P-2$, its positive inverse at zero is
\begin{equation*}
 Tf=(-A_{22})^{-1}f=R_2(f+\overline f),\qquad
 \overline{Tf}=\overline f,\qquad
 \norm{Tf}_\infty\le K\norm f_{L^1},\quad
 K=\frac2{1-e^{-2}}<4.
\end{equation*}
Indeed the mean identity follows by integration, substitution proves
$(-A_{22})Tf=f$, and the same calculation gives injectivity. These
identities also hold for signed continuous inputs.

\emph{Absence of a nonnegative eigenvector.}
Suppose $\mathcal A(u,v)=\lambda(u,v)$ with a nonzero nonnegative
domain vector, and let $m_u=\overline u$, $m_v=\overline v$.
If $m_u=0$, then $u=0$ and the first equation gives $v=0$.
Thus $m_u>0$. At $x=0$ the first equation is
$\lambda u(0)=m_u+v(0)>0$; it forces $u(0)>0$ and a real $\lambda>0$,
even if a complex eigenvalue was initially allowed. The second equation
then gives
\[
 m_v=\frac{m_u}{\lambda+1}\le m_u,\qquad
 v=R_{\lambda+2}(u+m_v),\qquad
 \norm v_\infty\le K m_u.
\]
On the other hand,
\begin{equation*}
 I_0:=\int_\T\frac{dx}{a(x)}
      =\frac2{16}\int_0^{1/2}(2x)^{-1/2}\,dx=\frac18.
\end{equation*}
The value at the single zero of $a$ is immaterial to this integral.
Solving the first eigenvalue equation and integrating yields
\[
 m_u=\int_\T\frac{m_u+v(x)}{\lambda+a(x)}\,dx
       \le(1+K)I_0m_u<\frac58m_u,
\]
a contradiction.

\emph{Growth bound.}
For every $\delta>0$ define
\[
 f_\delta=(\delta+a)^{-1},\qquad
 v_\delta=Tf_\delta,\qquad q_\delta=(f_\delta,v_\delta).
\]
This is a positive order unit in the actual mixed domain. Set
$I_\delta=\overline{f_\delta}\le I_0$. The second component of
$\mathcal A q_\delta$ is zero. The first obeys
\[
 (\mathcal A q_\delta)_1
 =\delta f_\delta-1+I_\delta+v_\delta
 \le\delta f_\delta-1+(1+K)I_\delta
 \le\delta f_\delta-\frac38.
\]
Thus $\mathcal A q_\delta\le\delta q_\delta$. The legitimate domain
identity and positivity give
\[
 e^{\delta t}q_\delta-S(t)q_\delta
 =\int_0^tS(t-s)e^{\delta s}
                   (\delta q_\delta-\mathcal A q_\delta)\,ds\ge0.
\]
Consequently $\norm{S(t)}\le M_\delta e^{\delta t}$ for a finite ratio
$M_\delta=\max q_\delta/\min q_\delta$. No uniform bound on this ratio
as $\delta\downarrow0$ is required, and $\omega_0(\mathcal A)\le0$.

\emph{Approximate spectrum and noncompactness.}
For integers $n\ge5$ let
\[
 \chi_n(x)=\max\{1-nr(x),0\},\qquad z_n=(\chi_n,T\chi_n).
\]
Then $\overline{\chi_n}=1/n$, $\norm{\chi_n}_\infty=1$ and
$\norm{T\chi_n}_\infty<4/n<1$, so $\norm{z_n}=1$.
Each $z_n$ is in $C\times C^1$, and
\[
 (\mathcal A z_n)_2=0,\qquad
 \norm{(\mathcal A z_n)_1}_\infty
       \le16\sqrt{2/n}+(1+K)/n\longrightarrow0.
\]
Hence zero is in the approximate spectrum. The Laplace-resolvent
inequality $s(\mathcal A)\le\omega_0(\mathcal A)$ now proves
\eqref{eq:zero-bound}. If a resolvent were compact, the bounded sequence
$\lambda z_n-\mathcal A z_n$ would make
$z_n=(\lambda-\mathcal A)^{-1}(\lambda z_n-\mathcal A z_n)$ have a
norm-convergent subsequence. Closedness of $\mathcal A$ would give a
nonzero nonnegative zero-eigenvector, already excluded. Thus the
resolvent is not compact.
\end{proof}

The counterexample identifies the loss of spatial control in the stationary
component. The reciprocal gap $1/a$ is integrable, so a concentrated spike
has small returned feedback through the moving component. It does not
contradict the all-nonzero-velocity theorem, whose resolvent estimate controls
every component derivative. Nor does it establish a counterexample for
Lipschitz potentials.

\begin{proposition}\label{prop:rankone}
Let $P$ be averaging on the unit circle, $k>0$, and
$a\ge0$ continuous with $\min a=0$. For $L=kP-a$ on $C(\T)$ set
$I(a)=\int_\T a(x)^{-1}dx\in(0,\infty]$. A strictly positive
eigenfunction exists if and only if $kI(a)>1$. In that case the
principal value is the unique $\lambda>0$ satisfying
\[
 k\int_\T\frac{dx}{\lambda+a(x)}=1;
\]
it is isolated and algebraically simple, its eigenfunction is
proportional to $(\lambda+a)^{-1}$, and it equals the growth bound.
If $kI(a)\le1$, then $s(L)=\omega_0(L)=0$ and there is no nonzero
nonnegative eigenfunction, including at equality.
For $a(x)=A(2\operatorname{dist}(x,\mathbb Z))^\zeta$, $A,\zeta>0$,
the threshold for $0<\zeta<1$ is $A(1-\zeta)<k$; for
$\zeta\ge1$ a positive principal function always exists.
\end{proposition}
\begin{proof}
The bounded operator generates a positive exponential: write
$L=-\|a\|_\infty\Id+[kP+(\|a\|_\infty-a)\Id]$ and use the
positive exponential series for the bracket.
At a zero of $a$, a nonzero nonnegative eigenfunction satisfies
$\lambda f=k\overline f>0$. Thus $\lambda>0$ and
$f=k\overline f/(\lambda+a)$. Monotone convergence identifies the
integral's limit at zero; dominated convergence gives its continuity
and differentiation for positive $\lambda$. It decreases strictly
from $kI(a)$ to zero, proving existence and nonexistence assertions.
For positive $\lambda$ the multiplication resolvent exists, and the
rank-one resolvent equation has scalar denominator
$1-k\int(\lambda+a)^{-1}$. Its derivative is positive. Thus the root
is isolated and simple; equivalently the same positive reciprocal
function is a dual density and excludes a Jordan vector by pairing.
Order-unit comparison gives both exponential norm bounds and hence
the asserted spectral and growth bound.

If $kI(a)\le1$, the order unit $q_\delta=(\delta+a)^{-1}$ obeys
$Lq_\delta\le\delta q_\delta$ for every $\delta>0$.
It follows that $\omega_0(L)\le0$. Choose continuous nonnegative
bumps of norm one supported on shrinking arcs about a minimum of $a$.
Their means and the supremum of $a$ on their supports tend to zero,
so $Lf_n\to0$ and zero is in the approximate spectrum. This proves
$s(L)=\omega_0(L)=0$. For the displayed power cusp,
$I(a)=1/[A(1-\zeta)]$ when $\zeta<1$, and $I(a)=\infty$ otherwise.
A nonnegative Lipschitz function vanishing somewhere satisfies
$a(x)\le C|x-x_0|$ nearby, so its reciprocal integral diverges.
This explains the regularity boundary in this particular rank-one
model; it is not a necessary-and-sufficient classification for general
coupled transport operators.
\end{proof}

\subsection{A resonant time-periodic obstruction}\label{sec:floquet}

The example below has a smooth strictly positive kernel, nonzero drift
in both species and coupling equal to one everywhere. Its potential is
H\"older continuous with a cusp; the obstruction is proved in this
regularity class.

On $\T=\R/\mathbb Z$ take the Gaussian--uniform kernel from the proof of
Theorem~\ref{th:zero}, so its convolution is $P$. Write
\[
 a(x)=16\sqrt{2\operatorname{dist}(x,\mathbb Z)},\qquad
 R=\begin{pmatrix}-1&1\\1&-1\end{pmatrix}.
\]
Consider the mild transport equation on $X=C(\T;\R^2)$
\begin{equation}\label{eq:time-model}
 u_t=\partial_xu+(P-\Id)u+B(t,x)u,\qquad
 B(t,x)=(1-a(x+t))\Id+R.
\end{equation}
Both velocities and both jump rates equal one. The matrix is continuous
and one-periodic in time and space, and both off-diagonal entries equal one.

\begin{theorem}[No positive Floquet function despite a uniform minorant]
\label{th:floquet}
The positive evolution family $U(t,s)$ of \eqref{eq:time-model} has a
one-period map $M=U(1,0)$ with $\spr(M)=1$ and no nonzero nonnegative
eigenvector. Consequently there are no real $\lambda$ and nonzero
nonnegative one-periodic mild function $p(t)$ such that $e^{\lambda t}p(t)$
solves \eqref{eq:time-model}. Nevertheless, for every $\tau>0$,
\begin{equation}\label{eq:time-minorant}
 (U(s+\tau,s)f)_i(x)\ge e^{-17\tau}\tau\sinh(\tau)
                     \sum_{j=1}^2\int_0^1f_j(y)\,dy
 \qquad(f\ge0)
\end{equation}
uniformly over $s,x,i$. The period map is not compact.
\end{theorem}
\begin{proof}
\emph{A scalar calculation.}
Let $L_0=P-a$ on $C(\T)$. It is bounded and generates a positive
semigroup, since $L_0=-16\Id+(P+16-a)$ with a positive term in parentheses.
There is no nonzero nonnegative eigenfunction. Indeed, if
$L_0f=\lambda f$, $f\ge0$, $f\ne0$, then $m_f=\int f>0$ and
$\lambda f(0)=m_f>0$, so $\lambda>0$, even if initially allowed complex.
The identity $f=m_f/(\lambda+a)$ would give
\[
 1=\int_0^1(\lambda+a(x))^{-1}\,dx
   <\int_0^1a(x)^{-1}\,dx=1/8,
\]
a contradiction. For $q_\delta=(\delta+a)^{-1}$,
\[
 L_0q_\delta=\delta q_\delta-1+\int q_\delta
       \le\delta q_\delta.
\]
This continuous order unit yields
$\norm{e^{tL_0}}\le(\max q_\delta/\min q_\delta)e^{\delta t}$
by positivity and variation of constants. The spikes
$\chi_n=(1-n\operatorname{dist}(x,\mathbb Z))_+$ have norm one,
mean $1/n$, and
$\norm{L_0\chi_n}\le n^{-1}+16\sqrt{2/n}\to0$.
Thus $s(L_0)=\omega_0(L_0)=0$.

\emph{The autonomous bounded operator.}
Set $A_0=L_0\Id+R$, $S_*(t)=e^{tA_0}$. It is positive because
\[
 A_0=-17\Id+Q,\qquad Q=(16-a)\Id+P\Id+N\ge0,
 \qquad N=\begin{pmatrix}0&1\\1&0\end{pmatrix}.
\]
The order unit $(q_\delta,q_\delta)$ gives
$\norm{S_*(t)}\le C_\delta e^{\delta t}$ for every $\delta>0$.
The vectors $z_n=(\chi_n,\chi_n)$ satisfy $\norm{z_n}=1$ and
$A_0z_n\to0$, so $s(A_0)=\omega_0(A_0)=0$.
Summing the two components of a nonnegative eigenvector of $A_0$ would
produce a nonnegative nonzero eigenfunction of $L_0$. Hence none exists.

\emph{The true mild evolution.}
For $(T(t)f)(x)=f(x+t)$, define
\begin{equation*}
 U(t,s)=T(t)S_*(t-s)T(-s).
\end{equation*}
This is a strongly continuous positive evolution family. To verify the
equation without assuming spatial differentiability of its solutions,
multiply the bounded autonomous Duhamel formula by $T(t)$. Since
$T(-r)a(\cdot+r)T(r)=a(\cdot)$ and $P,R$ commute with translations,
one obtains
\[
 u(t)=T(t-s)u_s+
       \int_s^tT(t-r)\{P-a(\cdot+r)\Id+R\}u(r)\,dr.
\]
The integrand is norm continuous. This is the mild transport equation
for \eqref{eq:time-model}; uniqueness follows from Gronwall's estimate.
No preservation of $C^1$ regularity across the cusp is needed. Since $T(1)=\Id$, $M=S_*(1)$.
The order-unit estimates at integer times give $\spr(M)\le e^\delta$
for every $\delta>0$. Conversely,
\[
 (M-\Id)z_n=\int_0^1S_*(r)A_0z_n\,dr\longrightarrow0.
\]
Therefore $1\in\sigma(M)$ and $\spr(M)=1$.

\emph{Exclusion of positive period eigenvectors.}
If $Mf=r f$ for $f\ge0$, $f\ne0$, positivity makes $r$ real
nonnegative; invertibility of $e^{A_0}$ makes $r>0$.
For $\lambda=\log r$, the Bochner integral
\[
 w=\int_0^1e^{-\lambda t}S_*(t)f\,dt
\]
is nonnegative and nonzero by continuity near zero. The boundedness of
$A_0$ justifies integration of the norm derivative, giving
\[
 (A_0-\lambda)w
   =[e^{-\lambda t}S_*(t)f]_{0}^{1}=e^{-\lambda}Mf-f=0.
\]
This contradicts the eigenvector exclusion for $A_0$. A nonnegative
periodic mild Floquet function would give exactly such an eigenvector
at time zero; if its value there vanished, uniqueness would make it zero.

\emph{The uniform minorant.}
Since $Q\ge P\Id+N$ and both are positive operators, induction yields
$Q^k\ge(P\Id+N)^k$. The operators $P\Id,N$ commute and $P^2=P$;
therefore
\[
 S_*(\tau)\ge e^{-17\tau}e^{\tau N}e^{\tau P\Id},\qquad
 e^{\tau P}=\Id+(e^\tau-1)P\ge\tau P,
 \qquad
 e^{\tau N}=\begin{pmatrix}\cosh\tau&\sinh\tau\\
                             \sinh\tau&\cosh\tau\end{pmatrix}.
\]
This proves \eqref{eq:time-minorant} for $S_*(\tau)$ and then for
$U(s+\tau,s)$, because translations preserve constants and integrals.
Finally, $z_n$ has no uniformly convergent subsequence: its pointwise
limit is discontinuous at zero. If $M$ were compact, a convergent
subsequence of $Mz_n$, together with $(M-\Id)z_n\to0$, would contradict
this fact.
\end{proof}

Each frozen operator in \eqref{eq:time-model} satisfies the autonomous
hypotheses of Theorems~\ref{th:principal} and \ref{th:mix}. The failure is
in the period evolution. Here the no-jump period map is
$V_0f(x)=e^{-a(x)}e^Rf(x)$. It is a contraction and
$(V_0-\Id)z_n\to0$, so $r(V_0)=1=r(M)$; the exact criterion
\eqref{tp:eq:spectral-gap} fails at equality.
Estimate~\eqref{eq:time-minorant} is unnormalized;
the proof of \eqref{eq:mix} first constructs a principal pair and then
normalizes to a Markov semigroup. These are different assertions for a
nonconservative evolution family.

\begin{corollary}\label{cor:one-active-no-gap}
Keep $a,R,P$ from \eqref{eq:time-model} and put $D=\diag(1,0)$.
The system
\[
 u_t=\partial_xu+D(P-\Id)u+[D-a(x+t)\Id+R]u
\]
has exactly one active jump component, positive constant off-diagonal
couplings and both characteristic trees. It has a uniform integral
minorant and a positive jump return, but no nonzero nonnegative
Floquet function. Its full and no-jump period radii both equal one.
\end{corollary}
\begin{proof}
Both assigned kernels are the smooth Gaussian--uniform density, the
velocities are one and the rates are $1,0$. The time-space coefficient
matrix is continuous and $1/2$-H\"older. At every frozen time the
full and no-jump generators have domain $C^1(\T;\R^2)$.
For $W=1$, both nontrivial edge scans are one. Take $h=17$ in
Lemma~\ref{oa:minorant-return}, with $m=2,n=1,\eta_a=1$.
It gives
\[
 U(s+3,s)f\ge e^{-51}\one\sum_i\int f_i,
 \qquad [U(s+2,s)-V(s+2,s)]\one\ge e^{-34}\one.
\]
These claims use no spectral separation.

In the moving coordinate the full bounded generator is
$A_1=\diag(P,0)-a\Id+R$. It generates a positive semigroup.
Suppose $A_1f=\lambda f$ with $f\ge0$, $f\ne0$.
Then $m_1=\int f_1>0$, since otherwise the first equation forces
$f_2=0$. At a zero of $a$, summing the equations gives
$\lambda(f_1+f_2)=m_1$, so $\lambda>0$. Eliminating the second
component gives
\[
 f_1=\frac{m_1}{2}\left\{
 \frac1{\lambda+a}+\frac1{\lambda+a+2}\right\}.
\]
After integration this would imply
$1<\frac12\int a^{-1}+1/4=5/16$, a contradiction
(the corresponding non-strict upper bound already suffices).
Thus $A_1$ has no such eigenfunction.

For $q_\delta=(\delta+a)^{-1}$, the vector $(q_\delta,q_\delta)$
is an order unit and obeys $A_1(q_\delta,q_\delta)\le
\delta(q_\delta,q_\delta)$, using $\int q_\delta\le1/8$.
Hence the growth bound is at most zero. The same vectors
$(\chi_n,\chi_n)$ as above have norm one and satisfy
$A_1(\chi_n,\chi_n)\to0$. Integration of the bounded semigroup
then gives $(e^{A_1}-\Id)(\chi_n,\chi_n)\to0$, so the period map
has radius one. If $e^{A_1}f=\rho f$ with $f\ge0$, $f\ne0$,
positivity and the minorant give $\rho>0$.
With $\lambda=\log\rho$, the nonzero nonnegative vector
$v=\int_0^1e^{-\lambda t}e^{tA_1}f\,dt$ satisfies
$(A_1-\lambda)v=e^{-\lambda}e^{A_1}f-f=0$, which is impossible.
The characteristic conjugacy from the previous proof transfers this
to absence of a periodic mild Floquet function.
Finally the no-jump map in that coordinate is multiplication by
$e^{-a(x)}e^R$. Its powers have norm one because $\min a=0$
and $e^{nR}\one=\one$. Thus its radius is one also.
\end{proof}

\subsection{Proof of the periodic realization proposition}\label{app:periodic-framework}
We verify the evolution-semigroup and positive-semigroup facts for the
transport realization of Section~\ref{sec:periodic}; see \cite[Propositions~3.1--3.2]{HR} and
\cite[Theorems~3.1 and~4.7]{GM}. The details make the characteristic domain
and both directions of the resolvent correspondence explicit.

\begin{proof}[Proof of Proposition~\ref{tp:th:general}]
\emph{The period boundary equation.}
The evolution extends to arbitrary ordered pairs of real times. Indeed,
the transport group is two-sided, and Picard--Dyson iteration for the
bounded continuous perturbation works in either direction, with terms
bounded by $C^n|t-s|^n/n!$. Uniqueness gives
$U(s,t)=U(t,s)^{-1}$, with common bounds on compact time blocks.
Backward operators are not assumed positive.

Write $D_0=C\partial_x-\partial_t$. Since
$(z-\mathscr G)w=F$ means $D_0w=zw-\mathcal B(t)w-F$,
the interaction variable $\mathcal T(-t)w(t)$ has derivative
$\mathcal T(-t)[\mathcal B(t)w(t)-zw(t)+F(t)]$ in $X$.
This fixes both the transport and forcing signs. For $F\in Y$,
variation of constants is therefore equivalent to
\begin{equation*}
 w(t)=e^{-zt}U(t,0)w(0)
       +\int_0^te^{-z(t-r)}U(t,r)F(r)\,dr.
\end{equation*}
The continuous directional derivative gives this formula for domain
functions. Conversely the characteristic integral equation gives that
derivative continuously; with the period endpoint condition it belongs
to $\mathcal D$. The directional derivative at the seam is
$zw-\mathcal B(t)w-F$, whose endpoint values agree because
$F,\mathcal B,w$ are periodic. Thus the periodic extension is in the
generator domain, rather than merely satisfying the endpoint equation. At $t=T$ the condition is
\begin{equation*}
 (\Id-e^{-zT}M)w(0)=H_zF,\qquad
 H_zF=\int_0^Te^{-z(T-r)}U(T,r)F(r)\,dr.
\end{equation*}
If $\Id-e^{-zT}M$ is invertible, these equations give a unique bounded
inverse of $z-\mathscr G$ with range in $\mathcal D$.

For the converse, $H_z:Y\to X$ has a bounded right inverse. Choose a
continuous scalar $\eta$ supported compactly in $(0,T)$, with
$\int_0^T\eta=1$. For $g\in X$ set
\[
 F_g(r)=\eta(r)e^{z(T-r)}U(r,0)M^{-1}g
\]
and extend periodically. It is continuous at the seam, its norm is
bounded by a constant times $\norm g$, and the evolution property gives
$H_zF_g=g$. If $z-\mathscr G$ is invertible, take its inverse on $F_g$
in the period-endpoint relation above to prove surjectivity of
$\Id-e^{-zT}M$. A nonzero vector in the latter kernel would give a
nonzero periodic kernel function of $z-\mathscr G$ by the homogeneous variation-of-constants formula above. Thus it is injective as well;
the bounded inverse theorem proves \eqref{tp:eq:spectrum-general}.

\emph{The growth bound.}
The positivity argument at the start of the proof of
Theorem~\ref{thm:compact-return} uses only positivity of $M$ and
$r(M)>0$, so $r(M)\in\sigma(M)$ here as well.
Identity~\eqref{tp:eq:spectrum-general} now places $\lambda_0$ in the
spectrum of $\mathscr G$ and excludes any spectral point with greater
real part.

For every $\varepsilon>0$, the spectral radius formula gives
$\norm{M^n}\le C_\varepsilon e^{(\lambda_0+\varepsilon)nT}$.
For $0\le s\le T$,
\[
 U(s+T,s)=U(s,0)M U(s,0)^{-1}.
\]
Both conjugating factors are uniformly bounded. Splitting off a residual
interval of length at most $T$ proves
$\norm{U(t,s)}\le C'_\varepsilon
 e^{(\lambda_0+\varepsilon)(t-s)}$ for all $t\ge s$.
At integer periods $U(nT,0)=M^n$, so no smaller uniform growth bound is
possible. By definition, the norm of $\mathcal E(\tau)$ is the supremum
of the corresponding evolution norms over the ending phase, proving
\eqref{tp:eq:growth-general}.

\emph{Upper variational formula and strict inverse.}
If $\mathscr Gq\le a q$ for an order unit, positive comparison bounds
the actual growth by $a$. Thus the infimum in
\eqref{tp:eq:upper-general} is at least $\lambda_0$.
For every real $\xi>\lambda_0$, its positive strong Laplace resolvent
gives $q_\xi=(\xi-\mathscr G)^{-1}\one$. A fixed short-time part of
the integral, using $\mathcal E(\tau)\one\ge e^{\alpha\tau}\one$,
proves $\min q_\xi>0$. Its quotient obeys
\[
 \frac{\mathscr Gq_\xi}{q_\xi}
       =\xi-\frac1{q_\xi}
       \le\xi-\frac1{\norm{q_\xi}}<\xi.
\]
Letting $\xi\downarrow\lambda_0$ proves the formula. This infimum need
not be attained.

If $\lambda_0<0$, the resolvent at zero supplies a strict positive
supersolution and a positive inverse. A strict supersolution conversely
has a strictly negative maximum quotient by compactness, proving
$\lambda_0<0$. If the zero inverse is positive, let
$q=(-\mathscr G)^{-1}\one\ge0$. It cannot vanish at a point of a
component: the continuous characteristic derivative is zero at that
minimum, while all convolution and off-diagonal terms there are
nonnegative, contradicting $(\mathscr Gq)_i=-1$.
Compactness makes $q$ an order unit and its quotient gives strict
stability. Finally the inverse proves the maximum principle:
$\mathscr Gf\ge0$ implies
$f=-(-\mathscr G)^{-1}\mathscr Gf\le0$.
\end{proof}

\section{Additional proofs for the nonlinear application}\label{app:applications}

\subsection{Proof of the nonlinear threshold theorem}\label{app:nonlinear-proof}
\paragraph{Well-posedness, comparison and scalar multiples.}
The square map is locally Lipschitz on $X$; transport Picard iteration
with a bounded clamp gives local mild solutions and continuous dependence
on the initial data by Gronwall. Comparison is obtained from a positive
Volterra equation. On any common cap $[0,M]$, choose $H$ so large that
$B_{ii}-d_i+H-2M\ge0$. Then the bounded map
$D_dKu+(B-D_d+H\Id)u-u\odot u$ is order preserving on that cap:
its convolution and off-diagonal terms are positive, and each diagonal
scalar map has nonnegative derivative. Using the base transport
$e^{-H(t-s)}\mathcal T(t-s)$, monotone Picard iteration or the Volterra
Gronwall estimate gives comparison for bounded solutions and sub/supersolutions.
The cap may be enlarged for a scalar multiple used below.
Let
$b_+=\max_{t,x,i}\sum_j B_{ij}(t,x)$.
For every
$M\ge\max\{\|u_0\|,\max(b_+,0)+1\}$,
the constant $M\one$ is a supersolution:
$D_d(K-\Id)\one=0$, $C\partial_x\one=0$, and
$B(M\one)-M^2\one\le M(b_+-M)\one\le0$.
Zero is a subsolution, so the solution remains in $[0,M\one]$
and continues globally. This also removes the temporary clamp.
The same comparison proves uniqueness for mild solutions.
The argument uses the transport integral formulation and requires no
separate spatial derivatives.

The function $p$ belongs to the continuous directional generator
domain and satisfies $p_t=A(t)p-\lambda p$ in that sense.
For a differentiable scalar $a(t)\ge0$, the residual of $a(t)p(t)$
in \eqref{nl:eq:nonlinear} is
\begin{equation*}
 (a'-\lambda a)p+a^2(p\odot p).
\end{equation*}
This identity can equally be integrated along each characteristic;
thus its use does not require $C^1$ spatial regularity of $p$.

When $\lambda>0$, every constant
$0<\varepsilon\le\lambda/p_+$ makes $\varepsilon p$
a lower solution. A large constant $M\one$ is the upper solution
already constructed. When $\lambda\le0$, choose
\[
 a'=\lambda a-p_-a^2,\qquad a(0)>0,\qquad u_0\le a(0)p(0).
\]
The residual is
$a^2p\odot(p-p_-\one)\ge0$, so $u(t)\le a(t)p(t)$.
For $\lambda=0$,
$a(t)=a(0)/(1+p_-a(0)t)$; for $\lambda<0$,
$a(t)\le a(0)e^{\lambda t}$. This proves the stated decay estimates.
If there is no $p$ but the actual uniform linear growth is negative,
the inequality $u(t)\le U(t,0)u_0$ proves decay directly.

\paragraph{Contraction in a positive order interval.}
Suppose $\lambda>0$.
Set $M_*=\max(b_+,0)+1$ and choose
\[
 0<\varepsilon_*<
 \min\{1,\lambda/(2p_+),M_* /(2p_+)\}.
\]
At phase zero let
$\mathcal I=[\varepsilon_*p(0),M_*\one]\subset X$.
The nonlinear period map preserves this interval, since the
lower profile is periodic and the upper constant is a supersolution.
Every trajectory in this interval remains above
$m_*=\varepsilon_*p_->0$ at all intermediate times.

On positive functions define the Thompson distance \cite{GQ}
\[
 d(f,g)=\max_{i,x}|\log f_i(x)-\log g_i(x)|.
\]
The interval is complete in this distance: its positive lower
bound and upper cap make this metric equivalent to the uniform
metric, and the interval is uniformly closed.
If two trajectories $v,w$ have the common floor $m>0$ and
$a_0=e^{d(v(0),w(0))}\ge1$, take
\[
 a'=-m a(a-1),\qquad a(0)=a_0.
\]
Since $A(t)$ is linear, the residual of $a(t)w(t)$ is
\[
 a'w+a(a-1)(w\odot w)
   =a(a-1)w\odot(w-m\one)\ge0 .
\]
Thus $v(t)\le a(t)w(t)$ by comparison; reversing the two
trajectories gives $w(t)\le a(t)v(t)$.
Moreover
$(\log a)'=-m(a-1)\le-m\log a$. Therefore
\begin{equation}\label{nl:eq:contraction}
 d(v(t),w(t))\le e^{-mt}d(v(0),w(0)).
\end{equation}
The period map on $\mathcal I$ is a strict contraction with factor
$e^{-m_*T}<1$. Banach's fixed-point theorem gives a continuous
positive periodic solution $u^*$.

Any two positive periodic solutions have positive joint floors
and finite caps by continuity on the compact period cell.
Applying \eqref{nl:eq:contraction} over $nT$ and using periodicity
forces their initial distance to be zero. This proves uniqueness
without assuming compactness or smoothing of the nonlinear map.

For any strictly positive initial function, reduce the lower
multiple of $p$ if necessary and enlarge the cap to contain both
its initial function and $u^*$. Both trajectories then have a
common fixed positive floor. Equation \eqref{nl:eq:contraction}
gives uniform exponential attraction. To convert distance to norm,
if both functions are at most $M$, then
$\|v-w\|\le M(e^{d(v,w)}-1)$.
Indeed $e^{d_0e^{-mt}}-1\le e^{d_0}d_0e^{-mt}$. Thus the constants
depend on this fixed cap and floor, as asserted.

\paragraph{Nonzero seeds and the entire-solution statement.}
On a cap $M$, $u_i^2\le M u_i$. Positive comparison with the linear
equation therefore gives
\[
 u(s+H)\ge e^{-MH}U(s+H,s)u(s).
\]
If \eqref{nl:eq:minorant} holds and $u(s)\ge0$ is nonzero, its
positive total integral makes $u(s+H)$ a strict order unit.
The preceding attraction argument applies at that phase.
Theorem~\ref{th:mix} supplies this minorant for autonomous coefficients
with every velocity nonzero and at least one active species.
Equation~\eqref{tp:eq:allminorant} supplies it for time-periodic coefficients
with every jump rate positive, and \eqref{oa:minorant} supplies it under
the two-tree test based at a designated active component.
Existence of a positive Floquet function
alone does not supply it when some components are inactive; see
Proposition~\ref{prop:smooth-lock}.
Without it, only strictly positive initial functions are asserted
to enter the positive contraction argument; pure translation with
zero initial values demonstrates why this restriction matters.

If $v(t)$ is an entire solution with
$0<m\le v_i(t,x)\le M<\infty$ for all $t,x,i$, compare it with
$u^*$ from time $t-S$ to any fixed time $t$.
The initial Thompson distance is at most a constant depending
only on the two global floors and caps. The contraction factor
is at most $e^{-m' S}$ for a fixed $m'>0$.
Letting $S\to\infty$ proves $v(t)=u^*(t)$.
This does not exclude pointwise-positive entire heteroclinics
whose joint infimum is zero.

If the coefficients are autonomous, every temporal translate of
$u^*$ is again a strictly positive $T$-periodic solution.
Uniqueness makes all such translates equal, hence $u^*$ is
time independent. Conversely any positive steady state is
periodic and falls in the same uniqueness class.

\subsection{The heterogeneous drift calculation}\label{app:drift-calculation}
\begin{proof}[Proof of Proposition~\ref{prop:drift-crossing}]
The kernel periodizes to $1/(2\pi)$ by the tiling calculation used in
Theorem~\ref{th:zero}; it has every moment and is smooth and strictly
positive. Both rates equal one and the two off-diagonal entries are
one. At $\eta=0$ the normalized principal function is $\one$, its
dual density has both components $1/(4\pi)$, and the principal value
is $b$. The fixed-domain resolvent proof of Theorem~\ref{th:perturb}
is analytic in these two finite-dimensional parameters. On a compact
velocity interval its Taylor remainders and their $c$ derivatives are
uniform. Translation by $\pi$ changes $\eta$ to $-\eta$, so the
principal value is even in $\eta$.

With the mean normalization of Theorem~\ref{th:principal}, write
\[
 \phi=\one+\eta\operatorname{Re}(w e^{ix})+O(\eta^2),\qquad
 \psi=\frac1{4\pi}
       [\one+\eta\operatorname{Re}(\overline w e^{ix})]+O(\eta^2).
\]
The conjugation in the second formula comes from reversing both
transport velocities in the independently constructed adjoint.
The first variation has zero value and its first Fourier mode solves
\[
 \begin{pmatrix}2-ic&-1\\-1&2+i\end{pmatrix}w=\binom10.
\]
The determinant is $(3+c)+i(2-2c)$ and
$w_1=[8+i(5c-1)]/(5c^2-2c+13)$. Pairing the second-order
equation with the constant dual density gives the coefficient
$\operatorname{Re}w_1/4$. This proves the first expansion; analytic
differentiation gives the second. Independently, the
Hellmann--Feynman formula \eqref{eq:HF} with
$\dot A=\operatorname{diag}(1,0)\partial_x$ gives at order $\eta^2$
$-\operatorname{Re}w_1\operatorname{Im}w_1/2
=-2(10c-2)/(5c^2-2c+13)^2$, the same coefficient. The constant-dual
pairing with every spatial derivative of a higher-order right mode
has zero mean, so no omitted term enters this coefficient.
At $c=1$ the coefficient after subtracting $3/32$ is $1/32$;
at $c=2$ it is $2/29-3/32=-23/928$. The uniform fourth-order
remainder proves the two strict signs and a crossing by continuity.
The derivative is strictly negative throughout $[1,2]$ for small
nonzero $\eta$, which gives uniqueness. The limiting equation for
the crossing is $15c^2-6c-25=0$; its positive root is the displayed
constant. The derivative there is nonzero, and the uniform remainder
in $s/\eta^2$ gives the stated $O(\eta^2)$ location error.

For the cubic assertion write $N=\left(\begin{smallmatrix}
2-ic&-1\\-1&2+i\end{smallmatrix}\right)$,
$D=\det N$ and $\Delta_c=|D|^2=5c^2-2c+13$.
The mean-normalized profile is $p=\one+\eta v+O(\eta^2)$,
$v=\operatorname{Re}(we^{ix})$, $w=N^{-1}(1,0)^T$.
Its first mean vanishes, so $a_*=1+O(\eta^2)$ and
$p^2-a_*p=\eta v+O(\eta^2)$. The stable inverse on the
complement of the simple principal projection is analytic and locally
uniform in $c$: explicitly it is the restriction of
$(A_{c,\eta}-s-\Pi)^{-1}(\Id-\Pi)$, whose inverse exists at
$\eta=0$ and persists by the fixed-domain resolvent identity.
Hence $r_*=-\eta\operatorname{Re}(N^{-1}we^{ix})+O(\eta^2)$.
Using $\ell(r_*)=0$ to subtract $a_*r_*$ in the cubic pairing gives
\[
 \mathfrak b
 =\frac{\eta^2}{2}\operatorname{Re}(\overline w^{\,T}N^{-1}w)
       +O(\eta^3),\qquad
 w=\frac{(2+i,1)^T}{D},\quad
 N^{-1}w=\frac{(4+4i,\,4+i(1-c))^T}{D^2}.
\]
The real pairing is $(2c^2+4c+58)/\Delta_c^2$.
Translation by $\pi$ preserves the mean normalization and conjugates
the profiles, duals and stable inverses at $\eta$ and $-\eta$;
the coefficient is therefore even. Analyticity improves the remainder
to uniform $O(\eta^4)$, proving \eqref{eq:drift-cubic} and its signs.

\end{proof}

\section*{Funding}
This work was supported by the National Natural Science Foundation of China
(Grant Nos.~12471176 and 12071491).
%The funding agency had no role in the study design, mathematical analysis, preparation of the manuscript, or decision to submit the article for publication.

\section*{Declaration of competing interests}
The author declares no known competing financial interests or personal relationships that could have influenced the work reported in this article.

\section*{Data availability}
No experimental or observational datasets were used in this theoretical study.
The proofs and explicit examples are included in the article.

\section*{Declaration of generative AI and AI-assisted technologies in the manuscript preparation process}
ChatGPT (OpenAI) was used during manuscript preparation for editorial revision, reference checking, and supplementary computational verification. The author reviewed and edited its output as necessary and takes full responsibility for the final manuscript.

\end{document}